\documentclass{amsart}
\usepackage{amsmath, amsthm, epsfig, eepic, epic}
\usepackage[knot,curve,dvips]{xy}

\newtheorem{thm}{Theorem}[section]
\newtheorem{cor}[thm]{Corollary}
\newtheorem{lem}[thm]{Lemma}
\newtheorem{prop}[thm]{Proposition}
\theoremstyle{definition}

\theoremstyle{remark}

\numberwithin{equation}{section}

\newenvironment{pf}{{\it Proof:}\quad}{\hfill$\Box$ \vskip 12pt}

\newcommand{\E}{E_K}
\newcommand{\G}{G(p,q)}
\newcommand{\Real}{\mathbb R}
\newcommand{\Ax}{A_{\kappa}}
\newcommand{\Aal}{A_{\al}}
\newcommand{\be}{\beta}
\newcommand{\al}{\alpha}
\newcommand{\ibe}{\beta^{-1}}
\newcommand{\ial}{\alpha^{-1}}
\newcommand{\als}{{\alpha^2}}
\newcommand{\bes}{\beta^2}
\newcommand{\ials}{\alpha^{-2}}
\newcommand{\ibes}{\beta^{-2}}
\newcommand{\alc}{{\alpha^3}}
\newcommand{\ialc}{\alpha^{-3}}
\newcommand{\ka}{\kappa^{p-18q}}
\newcommand{\imu}{\mu^{-1}}
\newcommand{\ma}{\mu\alpha^2\mu}
\newcommand{\ak}{\alpha^3\kappa^{p-18q}\alpha^3}
\begin{document}

\title{(-2,3,7)-pretzel knot and Reebless foliation}
\author{Jinha Jun}
\address{Department of Mathematics \\ Seoul National University \\ Seoul 151-747\\ Korea}%
\email{jhjun@math.snu.ac.kr}%
\thanks{The  author was partially supported by BK21}

\dedicatory{communicated by Rachel Roberts}%

% ----------------------------------------------------------------
\begin{abstract}
If $p/q > 18$, $p$ is odd, and $p/q\ne37/2$, $(p,q)$-Dehn surgery
for the (-2,3,7)-pretzel knot produces a 3-manifold without
Reebless foliation.
\end{abstract}
\maketitle
% ----------------------------------------------------------------
\section{introduction}

Every closed orientable 3-manifold admits a foliation with Reeb
components \cite{Ro}. On the contrary, Reebless foliation
$\mathcal{F}$ reflects the topological information of the ambient
manifold $M \supset \mathcal{F}$. Novikov\cite{No} showed that
leaves of $\mathcal{F}$ are $\pi_1$-injective and $\pi_2(M)=0$.
Rosenberg\cite{Ros} showed $M$ is irreducible or $M \approx S^2
\times S^1$. It follows that $\mathcal{F}$ lifts to
$\widetilde{\mathcal{F}}$ which has planar leaves in the universal
cover $\widetilde{M}$. Palmeira\cite{Pa} proved that any simply
connected $(n+1)$-manifold, $n\ge 2$, admitting a smooth
foliations by planar leaves with codim=1 is diffeomorphic to
$\Real^{n+1}$. It follows that the universal cover of $M$ is
homeomorphic to $\Real^3$(see also \cite{CC} for the proof).
Especially, $M$ is irreducible and $\pi_1(M)$ is infinite.

\begin{figure}[h]
\hspace{-1cm}
\begin{xy} /r1pc/:
0,{\vover,\vover-}
,+(4,2.5),{\vunder\vtwist\vunder-}
,+(4,5),{\vunder\vtwist\vtwist\vtwist\vtwist\vtwist\vunder-}
,(0,3.5),{\hcap[-3.5]}
,(9,3.5),{\hloop}
,(0,3.5);(9,3.5)**@{-}
%,(0,3.5);(9,3.5)**\crv{}
,(0,-2),{\hcap[-3.5]}
,(9,-4.5),{\hloop}
,(0,-5.5);(9,-5.5)**@{-}
%,(0,0);(9,2.5)**\crv{(0,5)&(9,3.5)}
,(1,0);(4,0.5)**\crv{(2,0)&(3,0.5)}
,(5,0.5);(8,2.5)**\crv{(6,0.5)&(7,2.5)}
%,(0,-2);(9,-4.5)**\crv{(0,-7)&(9,-5.5)}
,(1,-2);(4,-2.5)**\crv{(2,-2)&(3,-2.5)}
,(5,-2.5);(8,-4.5)**\crv{(6,-2.5)&(7,-4.5)}
\end{xy}
\caption{(-2,3,7)-pretzel knot}\label{pretzel}
\end{figure}
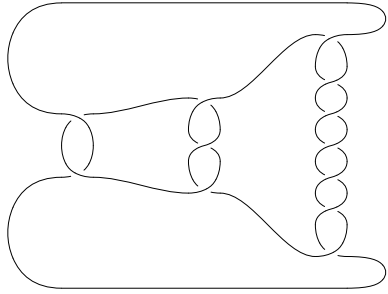

 Let $E_K$ be the (-2,3,7)-pretzel knot complement (Fig \ref{pretzel}). And let
$E_K(p/q)$ denote the 3-manifold obtained by $(p,q)$-Dehn surgery
along the (-2,3,7)-pretzel knot. It is known that there is no
closed essential surface in $\E$ and the boundary slopes are 0/1,
16/1, 37/2, and 20/1 \cite{HO}\cite{Oe}. Furthermore, $\E(16/1),
\E(37/2)$, and $\E(20/1)$ are toroidal \cite{HO}. $\E$ admits
(finite) cyclic surgery along $18/1$ and $19/1$\cite{FS}. And
Blieler and Hodgson\cite{BH} showed $\E(17/1)$ is a Seifert
fibered space with finite fundamental group.  In particular,
$\E(17/1), \E(18/1)$, and $\E(19/1)$ have no Reebless foliation by
virtue of Rosenberg's Theorem. $\E$ is fibered over the circle
with genus 5 surface whose monodromy is pseudo-Anosov and hence
hyperbolic. The suspension of the stable laminations gives an
essential lamination $\mathcal{L}$ in $\E$ with degeneracy slope
=1(18/1) \cite{Ga}. $\mathcal{L}$ remains essential in $\E(p/q)$
if $|p-18q|>1$ \cite{GO}. If $p$ is even, $\mathcal{L}$ extends to
a taut foliation in $\E(p/q)$ by filling  complementary regions
with a bundle of monkey saddle except $p/q =18/1$.

Using the technique in \cite{Ra} and \cite{Li}, one can prove that
$\E(p/q)$ has a Reebless foliation if $p/q \in (-\infty, 9)$. This
is done by attaching product disks to the fibers.

\begin{thm}[\textbf{Main Theorem}]
If $p/q > 18$, $p$ is odd, and $p/q \ne 37/2$, then $E_K(p/q)$
does not admit a Reebless foliation.
\end{thm}

The quotient space $\widetilde{M}/\widetilde{\mathcal{F}}$ is
called the leaf space. An open transversal to leaves gives an
1-manifold structure if $\mathcal{F}$ is a Reebless foliation. The
leaf space is a non-Hausdorff simply connected 1-manifold. There
is a natural action of $\pi_1(M)$ on the leaf space induced from
the action on $\widetilde{M}$. And this action has no global fixed
point (see \cite{Pa}). We will prove the Main theorem by showing
there is no nontrivial $\pi_1$-action on any leaf space.
 Our technique are much the same as in \cite{RSS}.

Calegari and Dunfield \cite{CD} notice that $\mathcal{F}$ gives
rise to a faithful  $\pi_1$-action on the a (universal) circle.
They showed there is no taut foliation in the Weeks manifold (the
closed hyperbolic 3-manifold with smallest known volume).

Our method is not applied to the case $p/q=37/2$. Indeed,
$\E(37/2)$ contains a Reebless foliation.  But this foliation is
not taut, because it has dead-end components. The following is
commented by Rachel Roberts.
\begin{lem}{\label{lem:37/2}}
$\E(37/2)$ does  contain  a Reebless foliation.
\end{lem}
\begin{pf}
Eudave-Mo\~{n}oz\cite{Eu} showed that $\E(37/2)$ is decomposed
along the incompressible torus $T^2$ into $\E(37/2) = X_L
\cup_{T^2} X_R$, where $X_L$ (respectively, $X_R$) is the
left-handed (respectively, right-handed) trefoil knot complement.

Since $X_L$ (respectively, $X_R$) is fibered, take the leaves of
the foliations which meet $\partial X_L =T^2 =\partial X_R$ in
simple closed curves of longitudinal slope and spiral them in a
neighborhood of the torus. By adding $T^2$ as a leaf, neither side
is a solid torus and so the resulting is Reebless.
\end{pf}

It is remarkable that any essential lamination in $\E(37/2)$
contains torus $T^2$ as a leaf \cite{BNR}.  The proof of Main
Theorem also can be used to show there is no transversely oriented
essential lamination except $p/q=37/2$. The following theorem
immediately follows from the results of \cite{RSS}.

\begin{thm}
If $p/q > 18$, $p$ is odd, and $p/q \ne 37/2$, $(p/q)$-Dehn
surgery for the (-2,3,7)-pretzel knot gives a 3-manifold without
transversely oriented essential lamination.
\end{thm}

If $M$ contains an essential lamination  with no isolated leaf, a
leaf space corresponds  to $\Real$-order tree \cite{GO}. In this
case, $\pi_1$ acts on $\Real$-order tree instead.

Since $\E(37/2)$ is Haken, it contains transversely oriented
essential lamination. In fact, there is the suspension of the
stable lamination in $\E(p/q)$ which remains essential when
$|p-18q|>1$. Main Theorem and the argument above imply the
following.

\begin{cor}
If $p/q>19$ and $p$ is odd, $\E(p/q)$ contains essential
lamination but does not admit any Reebless foliation.
\end{cor}

This paper is organized as follows. In Section \ref{sec:group}, we
discuss some basic properties of $\pi_1(\E(p/q)).$ Section
\ref{sec:action} gives an outline of theory of group actions on
(non-Hausdorff) simply connected 1-manifold. In Section
\ref{sec:real}, we will prove the nonexistence of $\Real$-covered
foliation in $\E(p/q)$. And the proof of the Main Theorem is
contained in Section \ref{sec:empty} and \ref{sec:nonempty}.

All results in this paper were obtained while the author was
visiting professor Rachel Roberts and professor John Shareshian in
the Washington University in 2002. This paper would not be
possible without their help. The author would like to express
thanks to them for their hospitality that makes visit to St. Louis
enjoyable and remarkable.

%------------------------------------------------------------------
\section{Fundamental Group}{\label{sec:group}}

This section contains useful properties  and a presentation of
(-2,3,7)-pretzel knot group. In later sections, we will analyze
the group actions on an orientable (non-Hausdorff) 1-manifold. The
following proposition implies the action can be restricted to the
orientation preserving one.
 Let $J$ be a knot in $S^3$ and $E_J$ be the exterior of $J$. Set
$G_0 = \pi_1(E_J(p/q))$.
\begin{prop}{\label{prop:index}}
If $p$ is odd, $G_0$ does not contain index 2 subgroup.
\end{prop}

\begin{pf}
 Suppose contrary that there is a subgroup $H$ with $[G_0:H]=2$.
 Since $G_0 /H\cong \mathbf{Z}
/2\mathbf{Z}$ is abelian, the commutator subgroup $[G_0,G_0] $ is
a subgroup of $H$. Note that  $G_0/[G_0,G_0]\cong H_1 (E_J(p/q)) =
\mathbf{Z}/p\mathbf{Z}$. Therefore we have a commutative diagram
below. Because $p$ is odd, we get a contradiction.

\begin{figure}[htb]
\setlength{\unitlength}{0.5cm}
\begin{picture}(10,4)
\put(0,3){$G_0$} \put(6,3){$G_0/[G_0,G_0]\cong
\mathbf{Z}/p\mathbf{Z}$} \put(2,0){$G_0/H \cong
\mathbf{Z}/2\mathbf{Z}$}

\put(2,3.2){\vector(1,0){3}} \put(1,2){\vector(1,-1){1}}
\put(7,2){\vector(-1,-1){1}}
\end{picture}
\end{figure}
\end{pf}

\begin{cor}{\label{cor:orient}}
Let $X$ be any oriented manifold and let
$$\Psi:G_0 \rightarrow Homeo(X)$$
be any homomorphism. Then $\Psi(G_0) \le Homeo^+(X)$.
\end{cor}

\begin{pf}
Suppose otherwise. Note that $[Homeo(X):Homeo^+(X)] \le 2$. Then
$\Psi^{-1}(Homeo^+(X))$ is an index 2 subgroup of $G_0$. By
Proposition \ref{prop:index}, it is impossible.
\end{pf}

Using the computer program SNAPPEA \cite{We}, we  can obtain a
presentation of the fundamental group of the knot (or link)
complement and the peripheral words using an ideal tetrahedra
decomposition. Denote (-2,3,7)-pretzel knot complement by $\E$.
The fundamental group of $\E$ and the meridian $m$ and longitude
$l$ are
$$\pi_1(E_K) = < a, b ~|~ a^2ba^2b^2a^{-1}b^2~>,$$
$$m=a^{-1}b^{-2}, \quad l=ab^{-1}a^2m^{-18}.$$
Of course, we have
\begin{equation}{\label{eqn:relation}}
a^2ba^2b^2a^{-1}b^2=1.
\end{equation}
 Let $G(p,q) :=\pi_1(E_K(p/q)) =<a,b ~|~
a^2ba^2b^2a^{-1}b^2, m^p l^q>.$
\begin{lem}[Lemma 3.4, \cite{RSS}]
There is some $k \in G(p,q)$ such that $m=k^q$ and $l=k^{-p}$.
\end{lem}
Now we have
\begin{equation}{\label{eqn:l}}
l=ab^{-1}a^2m^{-18} \Longleftrightarrow  k^{p-18q}= a^{-2}ba^{-1}.
\end{equation}

The following relation plays a central role in our proof.
\begin{equation}{\label{eqn:main}}
a^3k^{p-18q}a^3=a^3(a^{-2}ba^{-1})a^3=aba^2=a^{-1}b^{-2}ab^{-2}=ma^2m.
\end{equation}

%-------------------------------------------------------------------

\section{Group action on the leaf space}{\label{sec:action}}
We begin with a short exposition of the theory of group actions on
non-Hausdorff simply connected 1-manifold, taken from \cite{RSS}.
Let $\mathcal{F}$ be a Reebless foliation in $M$. Then $
\mathcal{F}$ can be lifted to $\widetilde{\mathcal{F}}$ in the
universal cover $\widetilde{M}$. The quotient space
$\mathcal{T}=\widetilde{M}/\widetilde{\mathcal{F}}$ is called the
\emph{leaf space}. The leaf space $\mathcal{T}$ is a simply
connected, 2nd countable 1-manifold \cite{CC}. But, in general, it
is not necessarily Hausdorff. Moreover there is an 1-1 correspond
between simply connected 1-manifolds and planar foliations in
$\Real^3$ up to conjugate by Palmeira \cite{Pa} (see also
\cite{CC} for details). Gabai and Kazez \cite{GK} extends this
relation to the essential laminations and $\Real$-order trees.

We recall here some terminology and definitions in \cite{RSS}.
Given $x,y \in \mathcal{T}$, we consider the \emph{geodesic spine}
$$[[x,y]]=\{z\in \mathcal{T}|x,y \mbox{ lie in distinct components of }
\mathcal{T}\setminus\{z\}\}\cup \{x,y\}$$
%\newline
from $x$ to $y$. $[[x,y]]$ is the union of a finite number of
disjoint (possibly, degenerate) closed intervals.
$$[[x,y]]=[x, y_1]\cup [x_2, y_2] \cup \dots \cup [x_n,y],$$
where $y_i$ is not separated from $x_{i+1}$. Set
$$d(x,y)=n-1.$$
Obviously, if $y\in [[x.z]]$ for some $x,y,z\in \mathcal{T}$
$$d(x,z)=d(x,y)+d(y,z).$$
Let us call a subset $X$ of $\mathcal{T}$ \emph{spine-connected}
if for all $x,y\in X$, $[[x,y]]\subset X$.

 Fix an orientation on $\mathcal{T}$. For
$x\in \mathcal{T}$, $\mathcal{T}\setminus \{x\}$ has exactly two
components since $\mathcal{T}$ is simply connected. If $U$ is a
connected Euclidean neighborhood of $x$, the two components of $U
\setminus\{x\}$ lie in distinct components of
$\mathcal{T}\setminus \{x\}$(Exercise C.1.4, \cite{CC}). Only one
component, say, $U^+$ is in the \emph{positive} direction of $x$.
Let $x^+$ be the component of $\mathcal{T}\setminus \{x\}$
containing $U^+$ and let $x^-$ be the component
$\mathcal{T}\setminus(x^+\cup \{x\})$.

Now we define a partial relation $\le$ on $\mathcal{T}$. For
$x,y\in \mathcal{T}$,
$$x \le y \Longleftrightarrow x^+ \supseteq y^+.$$
It follows that every map in $Homeo^+(\mathcal{T})$ preserves this
order $\le$.

Define a relation $\sim$ on $\mathcal{T}$ by $x\sim y$ if and only
if  $x $ and $y $ are not separated in $\mathcal{T}.$

Set
$$ [x]=\{y\in \mathcal{T} | y\sim x\}.$$
If $x\sim y$, let $\mathcal{T}_{\{x,y\}}$ denote the submanifold
defined as follows:
\begin{itemize}
\item if $x\in y^+$ (equivalently, if $y\in x^+$), set $\mathcal{T}_{\{x,y\}}=\bigcap_{z\sim x\mbox{ and } z\sim y}
z^+$, and
\item if $x\in y^-$ (equivalently, if $y\in x^-$), set $\mathcal{T}_{\{x,y\}}=\bigcap_{z\sim x\mbox{ and } z\sim y}
z^-$.
\end{itemize}

 The relation $\sim$ is
reflexive and symmetric, but not necessarily transitive. However,
by Denjoy blowing up, we can modify $\mathcal{F}$ so that $\sim$
is an equivalent relation (see Appendix in \cite{RSS}). In what
follows we shall assume $\sim$ is an equivalence relation.  Define
the \emph{Hausdorff tree} $\mathcal{T}_H=\mathcal{T}/\sim$.

If $X,Y$ are disjoint, nonempty, spine-connected subsets of
$\mathcal{T}$, the \emph{bridge} from $X$ to $Y$ is the
intersection of all paths in $\mathcal{T}$ with one end point in
$X$ and the other in $Y$. Similarly, we can define the bridge in
$\mathcal{T}_H$.

For any group $G$ acting on $\mathcal{T}$, if  $g\in G$, denote
$$Fix(g) = \{x\in \mathcal{T} | xg=x\}$$
and
$$Nonsep(g)=\{x\in \mathcal{T}| xg\sim x\}.$$
We say the action is \emph{trivial} or \emph{has a global fixed
point} if there is some $x \in \mathcal{T}$ such that $x\sim xg$
for all $g\in G$.

Define the \emph{characteristic set} associated to $g$ by
$$C_g=\{x\in \mathcal{T}|d(x,xg) \mbox{ is even }\}.$$

\begin{lem}[Lemma 4.7, \cite{RSS}]\label{lem:4.6}
Let $x\in \mathcal{T}$. Then $x\in C_g$ if and only if $x$ and
$xg$ are comparable with respect to the partial order $\le$.
\end{lem}

\begin{prop}[Proposition 4.8, \cite{RSS}]{\label{prop:4.7}}
Suppose $Nonsep(g)=\emptyset$. Then $C_g \ne \emptyset$ and for
any $x\in C_g$,
$$C_g=\bigcup_{n\in \mathbb{Z}} [[xg^n, xg^{n+1}]].$$
\end{prop}

When $Nonsep(g) =\emptyset$, $A_g :=C_g$ is called an \emph{axis}
for $g$. From  Proposition \ref{prop:4.7}, in $\mathcal{T}$, $A_g
\approx \Real$ or $A_g \approx \cup_{-\infty}^{\infty} [x_i,
y_i]$, where $[x_i, y_i]$ is homeomorphic to a closed interval in
$\Real$, $[x_i, y_i] \cap [x_{i+1}, y_{i+1}] = \emptyset$ when
$i\ne j$, $x_i\ne y_i$ and $y_i \sim x_{i+1}$ for all $i,j$. In
each case, the action of $g$ on $A_g$ is conjugate to an action by
translations. In $\mathcal{T}_H$, the image of $A_g$ is
homeomorphic to $\Real$.

Suppose $Y$ is a $g$-invariant embedded copy of $\Real$ in
$\mathcal{T}$ on which $g$ acts freely. Then we call $Y$ a
\emph{local axis} for $g$.  Now suppose that $Nonsep(g) \ne
\emptyset$ and let $T_i$ for some $i\in \mathcal{I}$, denote the
path components of $\mathcal{T}\setminus Nonsep(g)$. Notice that
$T_i g=T_j$ for some $j\in\mathcal{I}$. Moreover, whenever $T_i
g=T_i$, $g$ acts freely on $T_i$, and hence this local action has
an axis $A^i_g \subset T_i$. One can check that such an $A^i_g
\approx \Real$ and hence is an example of a local axis for $g$.
\begin{lem}[Lemma 4.10, \cite{RSS}]{\label{lem:lem4.9}}
Suppose $Nonsep(g) \ne \emptyset$. Then
$$C_g = Fix(g) \cup \{x\in \mathcal{T}| x \mbox{ lies on a local axis for }
 g\}.$$
\end{lem}
\begin{lem}[Corollary 4.12, \cite{RSS}]{\label{lem:cor4.11}}
If there is some $x \in \mathcal{T}$ such that $d(x, xg) \ne 0$ is
even, then $Nonsep(g) \ne \emptyset$.
\end{lem}

\begin{lem}[Corollary 4.13, \cite{RSS}]{\label{lem:cor4.12}}
Let $g\in G$. Then both $C_g$ and $C_g \cup Nonsep(g)$ are
spine-connected.
\end{lem}

Sometimes it is useful to consider an object obtained by adding
one points $\hat{x}$, called an \emph{ideal point} of
$\mathcal{T}$, to $\mathcal{T}$ for each $\sim$-equivalence class
$[x]$ in $\mathcal{T}$ which contains more than one point. This
object, denoted by $\widehat{\mathcal{T}}$, is called the
\emph{completion} of $\mathcal{T}$. We say that an ideal point
$\hat{x}$ is a \emph{source} if whenever $y,z$ are distinct
elements of $[x]$ we have $y\in z^-$ and we say that $\hat{x}$ is
a \emph{sink} if whenever $y,z$ are distinct elements of $[x]$ we
have $y\in z^+$. Note that every ideal point $\hat{x}$ is either a
source or a sink. The action of any subgroup of
$Homeo(\mathcal{T})$ extends to an action on
$\widehat{\mathcal{T}}$ in the obvious way, that is, we set
$\hat{x}g=\hat{y}$ if $[x]g=[y]$. We want to extend our partial
order on $\mathcal{T}$ to $\widehat{\mathcal{T}}$ so that group
actions on $\widehat{\mathcal{T}}$ obtained from orientation
preserving actions on $\mathcal{T}$ preserve this extended partial
order. For an ideal point $\hat{x}$, we define
$$\hat{x}^+=
\begin{cases}
\bigcup_{y\in [x]} (\{y\}\cup y^+), &\hat{x} \text{ a source, }\\
\bigcap_{y\in [x]} y^+, &\hat{x} \mbox{ a sink, }
\end{cases}$$
and set
$$\hat{x}^-=\mathcal{T}\setminus \hat{x}^+.$$
Note that $\hat{x}^+, \hat{x}^- \subset \mathcal{T}$. If $h\in
Homeo^+(\mathcal{T})$, for $x,y\in \widehat{\mathcal{T}}$, we have
$x^+ \subset y^+$ if and only if $(xh)^+ \subset (yh)^+$. So we
extend the partial order $\le$ in
$\mathcal{T}$ to $\widehat{\mathcal{T}}$. %Moreover orientation preserving map of $T$
%preserves the orientation on $\widehat{T}$.

Whenever possible, we will use $\mathcal{T}_H$ instead of
$\mathcal{T}$ to avoid tedious arguments when we deal with
non-Hausdorff points and to use the simply connectedness.
\begin{lem}[Lemma 5.6, \cite{RSS}]{\label{lem5.6}}
Any nontrivial action of $G$ on $\mathcal{T}$  canonically induces
a nontrivial action of $G$ on $\mathcal{T}_H$.
\end{lem}

 The action of $\pi_1(M)$ on $\widetilde{M}$ induces a right
action by homeomorphisms on $\mathcal{T}$. That is, there is a
homomorphism
$$\Phi: \pi_1(M) \rightarrow Homeo(\mathcal{T}).$$

By Corollary \ref{cor:orient}, we can assume $\Phi : G(p,q)
\rightarrow Homeo^+ (\mathcal{T})$.

We set
\begin{align*}
\begin{split}
&\al = \Phi(a)\\
&\be=\Phi(b)\\
&\mu=\Phi(m)\\
&\kappa=\Phi(k).
\end{split}
\end{align*}
Thus we have
\begin{alignat}{2}
\als\be\als\bes\ial\bes&=1{\label{eqn:rel}}&\quad \text{by (\ref{eqn:relation})},\\
\kappa^{p-18q}&=\ials\be\ial{\label{ka}}&\quad\text{by (\ref{eqn:l})},\\
\ak &=\ma{\label{eqn:word}}&\quad\text{by (\ref{eqn:main})}.
%\end{equation}
%$\because
%\ak=\alc\ials\be\ial\alc=\al\be\als=\ial\ibes\al\ibes=\mu\als\mu.$
\end{alignat}

%Denote $G:=\Phi(G(p,q))$ from now on.
By   \cite[Theorem 7.9]{CD}, we can assume $\Phi$ is injective,
that is, $\pi_1$ acts faithfully.
We will abuse the notation $G(p,q)$ for the image of $\Phi$. %We
%will say
When there is no ambiguity, we will simply say that $G(p,q)$,
instead of $\Phi(G(p,q))$, acts on a leaf space $\mathcal{T}$.
%-------------------------------------------------------------------
\section{$\Real$-covered foliation}{\label{sec:real}}

In this section, we will prove nonexistence of $\Real$-covered
foliation in $\E(p/q)$ for $p/q\ge10$.

The following lemma will be used in several times in the proof of
Main Theorem. If one wish to analyze other 3-manifold group
following \cite{RSS}, this lemma seems to be a criterion in
choosing a presentation of a group.
\begin{lem}{\label{poset}}
Let $G(p,q)$ act on a partially ordered set $P$. Suppose that
$G(p,q)$ preserves order. If some $x \in P$ satisfies either of
the conditions
\begin{enumerate}
\item $x\kappa = x$ and $x, x\al$ are related in $P$, or
\item $x\al = x$ and $x,x\kappa$ are related in $P$
\end{enumerate}
then $x$ is fixed by every $g \in G(p,q)$.
\end{lem}

\begin{pf}
Suppose $x\kappa = x$. Then we can assume $x < x\al$ because
$G(p,q) = <\kappa, \al>$. Then
\begin{equation}
x=x\mu=x\ial\ibes < x\ibes.{\label{eqn:1}}
\end{equation}
 But $x=x\ka = x\ials\be\ial <
x\be\ial$. This implies
$$ x\al < x\be.$$
Since $x < x\al$, we have $x < x\be < x\bes$. Contradiction to
(\ref{eqn:1}).

Similarly, if $x\al =x$ we may assume $x < x\kappa$. Then
\begin{equation}
x<x\mu=x\ial\ibes=x\ibes.{\label{eqn:2}}
\end{equation}
 But $x<x\ka =x\ials\be\ial =x\be\ial$.
This shows
$$x < x\be.$$
Hence $x < x\bes$. Contradiction to (\ref{eqn:2}).
\end{pf}
\begin{lem}\label{lem:beta}
If $q > 0$ and $x\kappa > x$ for all $x \in \Real $ then $x > x\beta $.
\end{lem}

\begin{pf}
Since $x\ial\ibes=x\mu=x\kappa^q > x$ for all $x \in \Real$, we
have
\begin{equation}
x\alpha^{-1} > x\beta^2.\label{eqn:3}
\end{equation}

Then
\begin{alignat*}{2}
x&=x\als\be\als\bes\ial\bes &\quad \text{by (\ref{eqn:rel})}\\
&<x\als\be\als\ial\ial\ial &\quad \text{by (\ref{eqn:3})}\\
&=x\als\be\ial
%\end{split}
\end{alignat*}
and $x\al\ibe < x\als$. By replacing $x\al$ with $x$ we get
\begin{eqnarray}
x\al &>& x\ibe \label{eqn:4}, \text{ and }\\
 x\alpha^2 &>& x\beta^{-2} \label{eqn:5}.
\end{eqnarray}

Thus
\begin{alignat*}{2}
x   &= x\alpha^2\beta\alpha^2\beta^2\alpha^{-1}\beta^2 &\quad\text{by (\ref{eqn:rel}}) \\
    &> x\ibes\be\ibes\bes\ial\bes &\quad\text{by (\ref{eqn:5}})\\
    &=x\ibe\ial\ibes &\\
    &>x\ibe\bes\bes &\quad\text{by (\ref{eqn:3}})\\
    &=x\beta^3&
\end{alignat*}
for all $x \in \Real$. This implies
$x > x \beta.$
\end{pf}

\begin{prop}{\label{prop:real}}
If $p/q \ge 10$ and $\phi : G(p,q) \longrightarrow Homeo^+(\Real)$ is
any homomorphism then there is some $ x \in \Real$ which is fixed by
every element of $\phi(G(p,q)).$
\end{prop}
\begin{pf}
By Lemma \ref{poset}, we may assume $x\kappa > x$ for all $x\in
\Real$.  Then $x > x\be$ for all $x\in \Real$ by Lemma
\ref{lem:beta}. %By Lemma \ref{poset}, $x\al > x$ or $x\al < x$ for
%all $x\in\Real$. If $x\al < x$,
Since
$$x=x\als\be\als\bes\ial\bes<x\als\als\ial=x\alc,$$
Lemma \ref{poset} implies that $x\al > x$ for all $x\in \Real$.

Since $x\mu > x$, $x > x\al\bes$ for all $x \in \Real$.
Therefore
\begin{equation}{\label{eqn:6}}
x=x\als\be\al(\al\bes)\ial\bes <x\als\be\al\ial\bes=x\als\beta^3.
\end{equation}

Now, we will prove the following.
\begin{equation}{\label{eqn:7}}
x\al > x\mu^2\mbox{  for all }x \in \Real.
\end{equation}

To see this, note that this is equivalent to $x\als\bes\al\bes >
x$.
%\begin{equation*}
\begin{alignat*}{2}
x(\als\bes)\al\bes &> x\ibe\al\bes &\mbox{ by (\ref{eqn:6}})\\
                   &> x\ibe \ibe \bes &\mbox{ by (\ref{eqn:4}})\\
                   &=x.
\end{alignat*}
%\end{equation*}
From (\ref{eqn:7}),
\begin{equation*}
\begin{split}
x\al > x\mu^2 &\Leftrightarrow x > x\mu^2\ial\\
              &\Leftrightarrow x\ibes > x\mu^2\ial\ibes=x\mu^3.
\end{split}
\end{equation*}
So we have
\begin{equation}{\label{eqn:8}}
x\ibes > x\mu^3.
\end{equation}

It follows that
\begin{alignat*}{2}
x\kappa^{18q-p} = x\al\ibe\als &> x\al\ibe\beta^{-3} &\mbox{ by (\ref{eqn:6})}\\
                               &> x\mu^2 \be^{-4} &\mbox{ by (\ref{eqn:7})}\\
                               &> x\mu^2 \mu^6 &\mbox{ by (\ref{eqn:8})}\\
                               &= x\mu^8\\
                               &= x\kappa^{8q}.
\end{alignat*}
%\end{equation*}
Since we assume $x\kappa > x$, $18q-p > 10q$. Hence $p/q < 10$.
Contradiction to the hypothesis.
\end{pf}
We suspect that the Proposition  is still true for $p/q\ge9$. On
the contrary, it is likely that the taut foliation for the
coefficients $p/q\in (-\infty, 9)$ are $\Real$-covered.

%-------------------------------------------------------------------
\section{$Nonsep(\kappa) = \emptyset$}{\label{sec:empty}}

In this section, we will show the Main Theorem when
$Nonsep(\kappa) =\emptyset$. From now, we will assume $p/q > 18
\Leftrightarrow p-18q > 0$, unless specified otherwise.

\begin{prop}\label{prop:1}
Suppose $Nonsep(\kappa) = \emptyset$. Then the action $G(p,q)$ on
$\mathcal{T}$ is trivial.
\end{prop}

\begin{pf}
Consider the action on $\mathcal{T}_H$. There are 3 cases for $\Ax
\cap \Ax\al$.
\begin{enumerate}
\item $\Ax \cap \Ax\al = \Ax$.
\item $\Ax \cap \Ax\al$ is a nonempty proper closed connected subset of $\Ax$.
\item $\Ax \cap \Ax\al=\emptyset$.
\end{enumerate}
For case (1), $\Ax \approx \Real$ is invariant under Im$\Phi$ and
hence there is a fixed point in $\mathcal{T}_H$ by Lemma
\ref{prop:real}. Thus, there is a global fixed point in
$\mathcal{T}$ by Lemma \ref{lem5.6}.

Case (2) will be proved in Lemma \ref{interval} and case (3) in
the Lemma below.
\end{pf}

\begin{lem}\label{lem:emptyset}
Suppose $Nonsep(\kappa) = \emptyset$. If $\Ax \cap \Ax\al =
\emptyset$, the action $\G$ on $\mathcal{T}$ is trivial.
\end{lem}

\begin{pf}
The relation (\ref{eqn:word}) applies to give
$$\Ax \als \mu =\Ax \mu \als \mu = \Ax \alc \ka \alc.$$
We will compare the bridges from $\Ax$ to $\Ax\als\mu$ and
$\Ax\ak$ to find a contradiction.

In $\mathcal{T}_H$, we define a total order on $\Ax$ by $x \preceq
x\kappa$ for all $x\in \Ax$.
Let $[r,s]$ be the bridge from $\Ax$ to $\Ax\al$ in $\mathcal{T}_H$. %Since
%there is the bridge $[r\al, s\al]$ from $\Ax\al$ to $\Ax\als$,
%Since $\Ax\al \cap \Ax\als=\emptyset$,  $\Ax = \Ax\als$ or $\Ax
%\cap \Ax\als=\emptyset$. If $\Ax = \Ax\als$, then $\Ax\be\ial =
%\Ax\als\ka = \Ax \Leftrightarrow \Ax\be = \Ax\al$. Then
%$\Ax\bes=\Ax\imu\ial=\Ax\ial=\Ax\al=\Ax\be$. So
%$\Ax\al=\Ax\be=\Ax$.
\begin{enumerate}
\item $\Ax \cap \Ax \als \ne \emptyset$.

Since the bridge from $\Ax\al$ to $\Ax\als$ is $[r\al, s\al]$, $s
=r\al$. Let $[x,y]=\Ax \cap \Ax \als$. Possibly,  $x$ or $y$
is not finite.
    \begin{enumerate}
    \item $\Ax \cap \Ax\alc =\emptyset$.

    Then the bridge from $\Ax$ to $\Ax\alc$ begins at $r$. Hence
    the bridge from $\Ax\alc$ to $\Ax\ak$ begins at $r\ka\alc$.
    From the Fig \ref{empty-0}, we see that the end point of
    the bridge from $\Ax\al$ to $\Ax\als$ is
    $$r\als=
    \begin{cases}
    r & \mbox{if } x\preceq r \preceq y,\\
    x & \mbox{if } r\preceq x,\\
    y & \mbox{if } y\preceq r,\\
    \end{cases}$$
    because the bridge from  $\Ax\al$ to $\Ax\als$ is $[r\al, r\als]$.
    On the other hand, the bridge from
    $\Ax\alc$ to $\Ax\als\mu$ begins at
    $$r\ka\alc=
    \begin{cases}
    r\al & \mbox{if } x\preceq r \preceq y,\\
    x\al & \mbox{if } r\preceq x,\\
    y\al & \mbox{if } y\preceq r.\\
    \end{cases}$$
    For all three cases, we have $r\ka\alc=r\alc \Leftrightarrow
    r\ka=r$. But $r \prec r\ka$. Contradiction.

\begin{figure}[htb]
\setlength{\unitlength}{0.5cm}
\begin{picture}(24,6)

\put(1,0){$x\preceq r \preceq y$}
\path(0,2)(5,2)\put(5,2){$\Ax$}
\path(0,4)(5,4)\put(5,4){$\Ax\al$}
\path(1,2)(0,1)\put(1,2.2){$x$}
\path(4,2)(5,1)\put(4,2.2){$y$}\put(5,1){$\Ax\als$}
\path(1,4)(0,5)\put(1,3.5){$x\al$}
\path(4,4)(5,5)\put(4,3.5){$y\al$}\put(5,5){$\Ax\alc$}

\put(10,0){$r\preceq x$}
\path(8,2)(13,2)\put(13,2){$\Ax$}
\path(8,4)(13,4)\put(13,4){$\Ax\al$}
\path(10,2)(9,1)\put(10,2.2){$x$}
\path(12,2)(13,1)\put(12,2.2){$y$}\put(13,1){$\Ax\als$}
\path(10,4)(9,5)\put(10,3.5){$x\al$}
\path(12,4)(13,5)\put(12,3.5){$y\al$}\put(13,5){$\Ax\alc$}

\put(17,0){$y\preceq r$}
\path(16,2)(21,2)\put(21,2){$\Ax$}
\path(16,4)(21,4)\put(21,4){$\Ax\al$}
\path(17,2)(16,1)\put(17,2.2){$x$}
\path(19,2)(20,1)\put(19,2.2){$y$}\put(20,1){$\Ax\als$}
\path(17,4)(16,5)\put(17,3.5){$x\al$}
\path(19,4)(20,5)\put(19,3.5){$y\al$}\put(20,5){$\Ax\alc$}

\Thicklines
\path(2.5,2)(2.5,4)
\put(2.5,1.5){$r$}
\put(2.5,4.2){$r\al$}
\path(9,2)(9,4)\put(9,1.5){$r$}\put(9,4.2){$r\al$}
\path(20,2)(20,4)\put(20,1.5){$r$}\put(20,4.2){$r\al$}
\end{picture}
\caption{$\Ax \cap \Ax\alc=\emptyset$}{\label{empty-0}}
\end{figure}

    \item $\Ax \cap \Ax\alc \ne\emptyset$.

    It follows that $\Ax\alc\cap\Ax\ak \ne \emptyset$. Equivalently,
    $\Ax\alc\cap\Ax\als\mu\ne\emptyset$. But
    $\Ax\als\cap\Ax\alc=\emptyset$. Let $[z,w]=\Ax\cap\Ax\alc$.
    Then $x\preceq y\prec z\preceq w$. Since the bridge from
    $\Ax\als$ to $\Ax\alc$ is $[r\als, r\alc]$, $y=r\als$ and
    $z=r\alc$. That is, $z=y\al \in \Ax\al$. In particular,
    $z\not\in \Ax$. Contradiction.
    \end{enumerate}

\item $\Ax \cap \Ax\als =\emptyset$

Then the bridge from $\Ax$ to $\Ax\als$ is $[r, s\al]$ (Fig \ref{empty-1}).
So the bridge from $\Ax$ to $\Ax \als \mu$ is $[r\mu,
s\al\mu]$(Fig \ref{empty-1}).

\begin{figure}[htb]
\setlength{\unitlength}{0.5cm}
\begin{picture}(20,6)

\path(0,1)(6,1) \put(6.2,1){$\Ax$} \path(0,3)(6,3)
\put(6.2,3){$\Ax\al$} \path(0,5)(6,5) \put(6.2,5){$\Ax\als$}

\put(2,0.5){$r$} \put(2,3.2){$s$} \put(3.8,2.5){$r\al$}
\put(3.8,5.2){$s\al$}

\path(12,1)(18,1) \put(18.2,1){$\Ax=\Ax\mu$} \path(12,3)(18,3)
\put(18.2,3){$\Ax\al\mu$} \path(12,5)(18,5)
\put(18.2,5){$\Ax\als\mu=\Ax\mu\als\mu$}

\put(13.8,0.5){$r\mu$} \put(13.8,3.2){$s\mu$}
\put(15.8,2.5){$r\al\mu$} \put(15.8,5.2){$s\al\mu$}
\put(9,3){\vector(1,0){1}} \put(9.2,3.2){$\mu$} \Thicklines
\path(2,1)(2,3) \path(4,3)(4,5) \path(2,3)(4,3) \path(14,1)(14,3)
\path(16,3)(16,5) \path(14,3)(16,3)

\end{picture}
\caption{bridge from $\Ax$ to $\Ax\ma$}{\label{empty-1}}
\end{figure}

    \begin{enumerate}
    \item $\Ax \cap \Ax\alc \ne \emptyset$.

    Since the bridge from $\Ax\als$ to $\Ax\alc$ is $[r\als,
    s\als]$ and $\Ax \cap \Ax\als =\emptyset$, $s\al=r\als \Leftrightarrow
    s=r\al$. See Figures \ref{empty-1} and \ref{empty-01}.
    Note that $[r,r\al]\cap [r\al, r\als]=\{r\al\}$.
    So $[r\al,r\als]\cap [r\als, r\alc]=\{r\als\}$. But the bridge
    $[r\als, r\alc]$ should contain $[r,r\als]$. Contradiction.

\begin{figure}[htb]
\setlength{\unitlength}{0.5cm}
\begin{picture}(10,4)

\path(1,1)(7,1)\put(7.5,1){$\Ax$}
\path(4,1)(3,0)
\path(5,1)(6,0)\put(6.5,0){$\Ax\alc$}
\path(1,2)(7,2)\put(7.5,2){$\Ax\al$}
\path(1,3)(7,3)\put(7.5,3){$\Ax\als$}

\Thicklines
\path(2,1)(2,3)
\put(2,0.5){$r$}
\put(2.1,2.2){$r\al$}
\put(2,3.2){$r\als$}

\end{picture}
\caption{$\Ax \cap \Ax\als =\emptyset$ and $\Ax \cap \Ax\alc \ne\emptyset$}{\label{empty-01}}
\end{figure}
    \item $\Ax \cap \Ax\alc = \emptyset$.

    The bridge from $\Ax$ to $\Ax \alc \ka \alc$ is $[r,s\als\ka\alc]$(Fig \ref{empty-2}).
\begin{figure}[htb]
\setlength{\unitlength}{0.5cm}
\begin{picture}(20,6)

\path(0,1)(6,1) \put(6.2,1){$\Ax\al^{-3}$} \path(0,3)(6,3)
\put(6.2,3){$\Ax=\Ax\ka$} \path(0,5)(6,5)
\put(6.2,5){$\Ax\alc\ka$}

\put(2,0.3){$r\al^{-3}$} \put(2,3.2){$s\ial$}
\put(3.3,2.3){$r\ka$} \put(3,5.2){$s\als\ka$}

\path(12,1)(18,1) \put(18.2,1){$\Ax$} \path(12,3)(18,3)
\put(18.2,3){$\Ax\alc$} \path(12,5)(18,5) \put(18.2,5){$\Ax\ak$}

\put(13.8,0.5){$r$} \put(13.8,3.2){$s\als$}
\put(15,2.2){$r\ka\alc$} \put(14,5.2){$s\als\ka\alc$}
\put(10,2){\vector(1,0){1}} \put(10.2,2.2){$\alc$}
 \Thicklines
\path(2,1)(2,3) \path(2,3)(4,3) \path(4,3)(4,5) \path(14,1)(14,3)
\path(16,3)(16,5) \path(14,3)(16,3)

\end{picture}
\caption{bridge from $\Ax$ to $\Ax\ak$}{\label{empty-2}}
\end{figure}

    Consequently,  $[r\mu, s\al\mu]=[r,s\als\ka\alc]$.
    But $r \ne r\mu=r\kappa^q$ in $\Ax \subset \mathcal{T}_H$.
    Contradiction.
    \end{enumerate}
\end{enumerate}
\end{pf}

Owing to the following lemma, we can rule out the case that
$\al$ acts on $\Ax$ with fixed points in $\mathcal{T}_H$ when
we prove Lemma \ref{lem:case2}
\begin{lem}
{\label{lem7.3}} Suppose that in $\mathcal{T}$, we have
$Nonsep(\al) \cap \Ax \ne \emptyset$. Then the action $\G$ on
$\mathcal{T}$ is trivial.
\end{lem}

\begin{pf}
If $x \in Fix(\al) \cap \Ax \ne \emptyset$, then $d(x, x\kappa)$
is necessarily even, and hence $x$ and $x\kappa$ are comparable
with respect to $\le$ on $\mathcal{T}$. Then Lemma \ref{poset}
applies. So
we may assume that $Fix(\al) \cap \Ax = \emptyset$, and choose $x \in Nonsep(\al) \cap \Ax$.\\
Consider $[x] \cap \Ax$. Either $[x] \cap \Ax = \{x\}$ or $[x]
\cap \Ax = \{x,y\}$ for some $y\ne x$. In the first case, Lemma
\ref{poset} applied to the ideal point determined by $[x]$ shows
that the action of $G(p,q)$ on $\mathcal{T}$ is trivial. Now we
assume the second case. Since $x \in \Ax$, we may assume $x\kappa
\le x$ (Fig \ref{fig:Ax}). Let $d(x,x\kappa) =2n> 0$.

\begin{figure}[hbt]
\setlength{\unitlength}{0.5cm}
\begin{picture}(15,2)
\put(2,1){\circle*{0.2}}
\put(0.7,0){$y\kappa^{-1}$}
\put(3,1){\circle*{0.2}}
\put(3,0){$x\kappa^{-1}$}
\put(7,1){\circle*{0.2}}
\put(6.5,0){$y$}
\put(8,1){\circle*{0.2}}
\put(8,0){$x$}
\put(12,1){\circle*{0.2}}
\put(11.5,0){$y\kappa$}
\put(13,1){\circle*{0.2}}
\put(13,0){$x\kappa$}

\put(0,1){\vector(1,0){1.9}}
\put(5,1){\vector(-1,0){1.9}}
\put(5,1){\vector(1,0){1.9}}
\put(10,1){\vector(-1,0){1.9}}
\put(10,1){\vector(1,0){1.9}}
\put(15,1){\vector(-1,0){1.9}}
\end{picture}
\caption{$Nonsep(\al) \cap \Ax \ne \emptyset $}{\label{fig:Ax}}
\end{figure}
There are 5 subcases:
\begin{enumerate}
\item $x=x\als$
\begin{align*}
\begin{split}
x > x\ka &= x \ials \be \ial = x\be \ial\\
&\Leftrightarrow x\al > x\be.
\end{split}
\end{align*}
But $x > x\mu = x\ial \ibes \Leftrightarrow x\bes > x\ial = x\al$.
Therefore,
\begin{align*}
\begin{split}
&x\bes > x\al  >x\be \\
&\Rightarrow x\be > x\\
&\Rightarrow x\al > x\be >x\\
&\Rightarrow x\al > x\\
&\Rightarrow x=x\als > x\al > x.
\end{split}
\end{align*}
Contradiction.
%\item $x=x\alc$ (and hence $x\ne x\als$)
\item  $y=x\als$\\
$x\be =x\als\ka\al=y\ka\al > x$ and $x\be \in (x\al)^-$ (Fig
\ref{fig:37/2}). Then $x\bes=x\imu\ial \in (y\ial)^-=(x\al)^-$ and
$x\bes> x\be$. We now see that $x\be \in [[x, x\bes]]$.
$d(x\be,x\bes) = d(x,x\be) = d(x,x\als\ka\al) = d(x\ial, y\ka) =
d(y,y\ka) = 2n(p-18q)$, because $x\ial \ne x$.

Therefore,
\begin{align*}
4n(p-18q) &= d(x, x\be) + d(x\be, x\bes)\\
        &= d(x, x\bes)\\
        &= d(x, x\imu\ial)\\
        &= d(x\al, x\imu)\\
        &= 2nq.
\end{align*}
That is, $p/q =37/2$.
However, we have assumed that $p/q \ne 37/2$.

%$E_K(37/2)$ is toroidal, it contains no Reebless foliation.

\begin{figure}[htb]
\setlength{\unitlength}{0.5cm}
\begin{picture}(17,8)
\put(0,1){$\Ax$}
\put(2,1){\circle*{0.2}}
\put(1,0){$x\imu$}
\put(6,1){\vector(-1,0){3.9}}
\put(6,1){\vector(1,0){3.9}}
\put(10,1){\circle*{0.2}}
\put(7.5,0){$x\als=y$}
\put(11,1){\circle*{0.2}}
\put(11,0){$x$}
\put(13,1){\vector(-1,0){1.9}}
\put(13,1){\vector(1,0){1.9}}
\put(15,1){\circle*{0.2}}
\put(14,0){$y\ka$}
\put(10.5,2){\circle*{0.2}}
\put(11,2){$x\al=y\ial$}
\put(10.5,3){\vector(0,-1){0.9}}
\put(10.5,3){\vector(0,1){0.9}}
\put(10.5,4){\circle*{0.2}}
\put(11,3.5){$x\be=y\ka\al$}
\put(10.5,5){\circle*{0.2}}
\put(10.5,6){\vector(0,-1){0.9}}
\put(10.5,6){\vector(0,1){0.9}}
\put(10.5,7){\circle*{0.2}}
\put(5,7){$x\imu\ial=x\bes$}
\end{picture}
\caption{$y=x\als$}{\label{fig:37/2}}
\end{figure}

\item $y\not=x\als$ and $x=x\alc$ (hence $x\ne x\als$)\\
 $d(x,x\alc\ka\alc) = d(x,x\ka)=2n(p-18q)$. On the other hand, $d(x,
x\mu\als\mu) = 4nq$ (Fig \ref{lem7.3-fig1}). Hence $2n(p-18q) =
4nq \Leftrightarrow p=20q$. But we assume $p$ is odd.
    %\end{enumerate}

%\item $x \ne x\als, x\alc$ and $y=x\als$ (and hence $y \ne x\alc$)
\item $x \ne x\als, x\alc$ and $y=x\alc$ (hence $y\ne x\als$)\\
$x\ak=y\ka\alc \in (x\alc)^-=y^-$, since $y\ka \in x^-$.  But
$x\ma \in y^+$ (Fig \ref{lem7.3-fig1}).
\item $x\ne x\als, x\alc$ and $y \ne x\als, x\alc$\\
$x\ma \in x^-$, but $x\ak \in x^+$ (Fig \ref{lem7.3-fig1},
\ref{lem7.3-fig2}).
\end{enumerate}
\end{pf}

\begin{figure}[htb]
\setlength{\unitlength}{0.5cm}
\begin{picture}(12,6)
\put(10,1){$\Ax$}
\put(1,1){\circle*{0.2}}
\put(0,1){\vector(1,0){0.9}}
\put(1,0){$y$}
\put(2,1){\circle*{0.2}}
\put(1.8,0){$x$}
\put(4,1){\vector(-1,0){1.9}}
\put(4,1){\vector(1,0){2.9}}
\put(7,1){\circle*{0.2}}
\put(6.5,0){$y\mu$}
\put(8,1){\circle*{0.2}}
\put(7.7,0){$x\mu$}
\put(1.5,1.5){\circle*{0.2}}
\put(1.8,1.5){$x\als$}
\put(1.5,3.5){\circle*{0.2}}
\put(1.5,4.5){\circle*{0.2}}
\put(1.8,4.5){$x\mu\als$}
\put(7.5,1.5){\circle*{0.2}}
\put(7.8,1.5){$x\als\mu$}
\put(7.5,3.5){\circle*{0.2}}
\put(7.5,4.5){\circle*{0.2}}
\put(7.8,4.5){$x\mu\als\mu$}

\path(1.5,1.5)(1.5,3.5)
\path(1.5,4.5)(1.5,6)
\path(7.5,1.5)(7.5,3.5)
\path(7.5,4.5)(7.5,6)
\put(9.5,1){\vector(-1,0){1.4}}

\end{picture}
\caption{$x\ma$}{\label{lem7.3-fig1}}
\end{figure}

\begin{figure}[htb]
\setlength{\unitlength}{0.5cm}
\begin{picture}(12,6)
\put(2,1){\circle*{0.2}}
\put(2,0){$y$}
\put(0,1){\vector(1,0){1.9}}
\put(3,1){\circle*{0.2}}
\put(3,0){$x$}
\put(6,1){\vector(-1,0){2.9}}
\put(6,1){\vector(1,0){1.9}}
\put(8,1){\circle*{0.2}}
\put(6,0){$y\ka$}
\put(9,1){\circle*{0.2}}
\put(8.7,0){$x\ka$}
\put(2.5,1.5){\circle*{0.2}}
\put(2.8,1.5){$x\alc$}
\put(2.5,3.5){\circle*{0.2}}
\put(2.5,4.5){\circle*{0.2}}
\put(8.5,1.5){\circle*{0.2}}
\put(8.8,1.5){$x\alc\ka$}
\put(3,4){\circle*{0.2}}
\put(3.2,4.2){$x\ak$}
\put(2.5,3.5){\vector(0,-1){1.9}}

\path(2.5,4.5)(2.5,6)
\path(3,4)(6,4)
\path(8.5,1.5)(8.5,4)
\put(11,1){\vector(-1,0){1.9}}

\put(11.2,0.8){$\Ax$}
\end{picture}
\caption{$x\ak$}{\label{lem7.3-fig2}}
\end{figure}

The following lemma, together with Lemma \ref{lem:emptyset}, will
complete the proof of Proposition \ref{prop:1}.
\begin{lem}\label{lem:case2}
{\label{interval}} Suppose that $\Ax \cap \Ax\al$ is  a nonempty
proper closed connected subset  of $\Ax$. Then the action $\G$ on
$\mathcal{T}$ is trivial.
\end{lem}

\begin{pf}
Suppose otherwise. In $\mathcal{T}_H$, $\Ax \approx \Real$. Let
$\preceq$ denote the total order on $\Ax$ specified by $x \preceq
x\kappa$ for all $x \in \Ax$. With respect to this order, let $r$
(respectively, $s$) denote the lower bound (respectively, upper
bound), if it exists, of $\Ax \cap \Ax\al$. Otherwise, set
$r=-\infty$ (respectively, $s=\infty$). Note that at least
one of $r$ and $s$ is finite because the intersection is a proper subset.\\
When $r \ne s$ let $\preceq_\al$ denote a total order on $\Ax\al$
such that $\preceq$ and $\preceq_\al$ agree on $\Ax \cap \Ax\al$.
Similarly, we define $\preceq_\als$ (respectively, $\preceq_\alc$)
on $\Ax \als$ (respectively, $\Ax\alc$) to agree with
$\preceq_\al$ (respectively, $\preceq_\als$) on $\Ax\al \cap
\Ax\als$ (respectively, $\Ax\als \cap \Ax\alc$).

Let $Y$ be an embedded copy of $\Real$ in $\mathcal{T}_H$ with a
total order. Then the homeomorphism
$$\al^r : Y \rightarrow Y\al^r$$ is order-preserving or
order-reversing if some total order is defined in $Y\al^r$.
 Suppose that $\al : (\Ax, \preceq) \rightarrow (\Ax\al,
\preceq_\al)$ is order-preserving. Then $r\al \preceq_\al s\al$ and
$r\ial \preceq s\ial$. Since $r\al, s\al \in \Ax\als$, $r\al
\preceq_\als s\al$. Therefore the map $\al : (\Ax\al, \preceq_\al)
\rightarrow (\Ax\als, \preceq_\als)$ is order preserving. In
particular, $r\als \preceq_\als s\als \Leftrightarrow r\als
\preceq_\alc s\als$. Since $(r\ial)\alc = r\als \preceq_\alc s\als
= (s\ial)\alc$, the map $\alc : (\Ax, \prec) \rightarrow (\Ax\alc,
\prec_\alc)$ is order preserving.

If $r$ and $s$ are finite and either
$$[r\al, s\al] \subset [r,s] \mbox{ or } [r,s] \subset [r\al,
s\al],$$ the intermediate value theorem implies $Nonsep(\al) \cap
\Ax \ne \emptyset$. This contradicts to Lemma \ref{lem7.3}.

Accordingly, we have 3 cases by symmetry (Fig \ref{order preserving}),\\
(I) $ r \prec_\al r\al \prec_\al s \prec_\al s\al$.\\
(II) $ r \prec_\al r\al = s \prec_\al s\al$.\\
(III) $   r \preceq_\al s \prec_\al r\al \preceq_\al s\al$ (if
$\preceq_\al$ is defined) $\Leftrightarrow$ $[r,s] \cap [r\al,
s\al]= \emptyset$ $\Leftrightarrow$ $\Ax \cap \Ax\als =
\emptyset$.

\begin{figure}[hbt]
\setlength{\unitlength}{0.5cm}
\begin{picture}(20,12)

\path(7,9)(12,9)
\put(5.5,6.5){$\Ax$}
\put(8,7){(II)}
\put(5.5,11.5){$\Ax\al$}
\path(7,9)(6,7)
\put(7.2,8.5){$r$}
\put(9.7,8.5){$r\al=s$}
\path(7,9)(6,11)
\path(12,9)(13,7)
\path(12,9)(13,11)

\path(2,3)(1,1)
\path(2,3)(1,5)
\path(2,3)(7,3)
\path(7,3)(8,1)
\path(7,3)(8,5)
\put(0.5,0.5){$\Ax$}
\put(3,0){(I)}
\put(0.5,5.5){$\Ax\al$}
\put(2.2,2.5){$r$}
\put(6.5,2.5){$s$}
\put(4.5,2.5){$r\al$}
\put(7.7,4.2){$s\al$}
\put(5,4.5){$\Ax\als$}

\path(12,3)(11,1)
\path(12,3)(11,5)
\path(12,3)(17,3)
\path(17,3)(18,1)
\path(17,3)(18,5)

\put(10.5,0.5){$\Ax$}
\put(13,0){(III)}
\put(10.5,5.5){$\Ax\al$}
\put(12.2,2.5){$r$}
\put(16.5,2.5){$s$}
\put(17.5,3.5){$r\al$}
\put(18,4.3){$s\al$}
\put(15,5){$\Ax\als$}
\Thicklines
\path(5,3)(5,4)
\path(5,3)(7,3)
\path(7,3)(7.5,4)
\path(7.5,4)(6.5,4.5)
\path(17.33,3.66)(16.33,4.66)
\path(17.66,4.33)(16.66,5.33)
\path(17.33,3.66)(17.66,4.33)

\path(12.5,10)(11,10.5)\put(13,10){$s\al$}
\put(10,11){$\Ax\als$}
\path(12.5,10)(12,9)
\path(12,9)(12.5,8)\put(13,8){$z$}
\path(12.5,8)(11,7.5)
\end{picture}
\caption{$\Ax, \Ax\al,\text{ and }\Ax\als$} {\label{order
preserving}}
\end{figure}
For (I), note that
$$r\al \prec_\als r\als \prec_\als s\al \prec_\als s\als. $$
\begin{figure}[htb]
\setlength{\unitlength}{0.5cm}
\begin{picture}(20,6)
\path(0,1)(17.5,1)
\put(18,0.8){$\Ax$}
\path(2,1)(2,5)
\put(1.5,5.2){$\Ax\als$}
\path(8,1)(8,5)
\put(7.5,5.2){$\Ax\als$}
\dottedline{.1}(15,1)(15,2)
%\dashline[+30]{3}(15,1)(15,2)
\put(17,3.7){$\Ax\alc\ka$}
\path(8,3)(6,3)
\put(4.3,3){$\Ax\al$}
\put(1.7,0){$r\al$}
\put(8,0){$s$}
\put(8.2,3){$s\al$}
\put(6.7,2){$r\als$}
\put(6.7,4){$s\als$}
\put(11.8,2){$r\als\ka$}
%\put(11.8,4){$s\als\ka$}
\put(14,0){$s\ka$}
\put(10,3.7){$\Ax\alc$}

\Thicklines
\path(8,2)(10,2)
\path(8,4)(10,4)
\path(8,2)(8,4)

\path(15,2)(17,2)
\path(15,4)(17,4)
\path(15,2)(15,4)

\end{picture}
\caption{$\Ax \cap \Ax\alc = \emptyset$}{\label{I-fig1}}
\end{figure}

If $\Ax \cap \Ax\alc = \emptyset$, then $s \prec_\als r\als$. Then
the bridge from $\Ax\alc$ to $\Ax\ak$ is a translation of the
bridge from $\Ax$ to $\Ax\alc\ka$ by $\alc$.
$$[s\ka, r\als\ka]\alc = [s\ka\alc, r\als\ka\alc]\text{(see Fig \ref{I-fig1})}.$$ On the other
hand, the bridge from $\Ax\alc$ to $\Ax \ma$ is
\begin{enumerate}
\item $[r\als, s]$ when $r\al\mu \prec s$,
\item $[r\als, r\al\mu]$ when $s \prec r\al\mu$,
\item $[r\als, x]$ when $s=r\al\mu$ for some $x$.
\end{enumerate}

See Fig \ref{I-fig2}.

\begin{figure}[hbt]
\setlength{\unitlength}{0.5cm}
\begin{picture}(20,6)

\path(2,2)(2,4)
\path(2,2)(4,2)
\path(2,4)(4,4)
\put(4,3.7){$\Ax\alc$}
\put(1,2){$r\als$}
\put(0,0.5){$r\al\mu$}
\put(2,0.5){$s$}
\dottedline{0.1}(2,1)(2,2)
\put(4.5,0.5){$s\mu$}

\put(0,5.2){$\Ax\als\mu$}

\path(8,2)(8,4)
\path(8,2)(10,2)
\path(8,4)(9,4)
\put(7,2){$r\als$}
\put(8,0.5){$s$}
\dottedline{0.1}(8,1)(8,2)
\put(10,0.5){$r\al\mu$}
\put(12,0.5){$s\mu$}
\put(8,4.2){$\Ax\alc$}
\path(0,1)(6,1)

\path(7,1)(13,1)

\path(14,1)(20,1)
\path(16,1)(15,5)
\path(15.5,3)(17,3)
\path(15.25,4)(17,4)
\put(14.3,3){$r\als$}
\put(17,3.5){$\Ax\alc$}
\put(15,0.5){$s=r\al\mu$}
\put(19,0.5){$s\mu$}
\put(16,1.8){$x$}

\Thicklines
\path(1,1)(0,5)
\path(5,1)(6,5)
\path(1,1)(5,1)
\path(11,1)(10,5)
\path(12,1)(13,5)
\put(10.5,5){$\Ax\als\mu$}
\path(11,1)(12,1)
\path(16,1)(15.75,2)
\path(19,1)(20,5)
\put(19,5.2){$\Ax\als\mu$}
\path(15.75,2)(15,2.2)
\path(16,1)(19,1)
\end{picture}
\caption{bridge from $\Ax\alc$ to $\Ax\ma$}
{\label{I-fig2}}
\end{figure}

Case (1) , $\begin{cases}
s\ka\alc = r\als,\\
r\als\ka\alc=s.
\end{cases}$ Hence both $r$ and $s$ are finite. By (\ref{ka}), the
second relation is equivalent to  $s=r\be\als$. Substituting this
in the first relation yields
\begin{align*}
\begin{split}
 r\als &= r\be\als\ka\alc =r\be\als(\ials\be\ial)\alc\\
 &\Leftrightarrow r = r\bes =r\imu\ial\\
 &\Leftrightarrow r\al=r\imu.
\end{split}
\end{align*}
But $r\imu \prec r \prec r\al$. Contradiction.

Case (2), $\begin{cases}
s\ka\alc = r\als,\\
r\als\ka\alc = r\al\mu.
\end{cases}$
Hence both $r$ and $s$ are finite. We have
\begin{align}\label{eqn:contradiction}
\begin{split}
r\al\mu&=r\als(\ials\be\ial)\alc=r\be\als\\
&\Leftrightarrow
r=r\be\als\imu\ial=r\be\als\bes=r\ials\ibes\al=r\ial\mu\al\\
&\Leftrightarrow r\ial = r\ial\mu.
\end{split}
\end{align}
But $r\ial \prec r\ial\mu.$ Contradiction.

Case (3),
$\begin{cases}
s\ka\alc=r\als,\\
r\als\ka\alc=x,\\
s=r\al\mu.
\end{cases}$ Then both $r$ and $s$ are finite. We have
\begin{align*}
\begin{split}
&r\al\mu\ka\alc = r\als\\
&\Leftrightarrow r\al\kappa^{p-17q} = r\ial.
\end{split}
\end{align*}
But $r\ial \prec r \prec r\al\kappa^{p-17q}.$ Contradiction.

Now we may assume $\Ax \cap \Ax\alc \ne \emptyset$ (Fig \ref{I-fig3}).
\begin{figure}[hbt]
\setlength{\unitlength}{0.5cm}
\begin{picture}(20,6)
\path(0,1)(15,1)
\put(15.5,0.7){$\Ax$}
\path(2,1)(2,5)
\put(1.8,0){$r\al$}
\put(1.5,5.3){$\Ax\als$}
\path(9,5)(9,1)
\put(8.5,5.3){$\Ax\al$}
\put(9,0){$s$}
\path(9,3)(12,3)
\put(12,2.8){$\Ax\als$}
\put(8,2.8){$s\al$}

\put(4.5,3.2){$\Ax\alc$}
\put(5.5,0){$r\als$}
\put(9.5,2.3){$s\als$}
\put(11,5.2){$\Ax\alc$}
\Thicklines
\path(5.5,3)(6,1)
\path(6,1)(9,1)
\path(9,1)(9,3)
\path(9,3)(10,3)
\path(10,3)(11,5)

\end{picture}
\caption{$\Ax \cap \Ax\alc  \ne \emptyset$}{\label{I-fig3}}
\end{figure}
Then $r\als\preceq_\als s$. Hence the intersection between $\Ax$
and $\Ax\alc\ka$ is $[r\als\ka, s\ka]$ in $\prec$-order in $\Ax$.
Since $\alc : (\Ax, \prec) \rightarrow (\Ax\alc, \prec_\alc)$ is
an order-preserving map, the intersection between $\Ax\alc$ and
$\Ax\alc\ka\alc$ is $[r\als\ka\alc, s\ka\alc]$ in
$\prec_\alc$-order in $\Ax\alc$. On the other hand, the
intersection between $\Ax\alc$ and $\Ax \als\mu$ is, in
$\prec_\alc$-order,
\begin{enumerate}
\item $[r\als, s]$ when $r\al\mu \prec r\als$,
\item $[r\al\mu,s]$ when $r\als\prec r\al\mu \prec s$,
\item $[x,s]$ when $r\als=r\al\mu$,
\item $[s,y]$ when $s=r\al\mu$.
\end{enumerate}
See Fig \ref{I-fig4}.

\begin{figure}[hbt]
\setlength{\unitlength}{0.5cm}
\begin{picture}(18,12)
\path(0,9)(7,9)
\put(2,7){$r\al\mu \prec r\als$}
\put(8,8,7){$\Ax$}
\path(1,9)(1,11)
\put(0.5,8.5){$r\al\mu$}
\put(0,11.2){$\Ax\als\mu$}
\path(6,9)(6,11)
\put(5.5,8.5){$s\mu$}
\path(10,9)(17,9)
\put(12,7){$r\als \prec r\al\mu \prec s$}
\path(13,9)(13,11)
\put(12.5,8.5){$r\al\mu$}

\path(16,9)(16,11)
\put(16,8.5){$s\mu$}
\put(16,11.2){$\Ax\als\mu$}
\path(0,3)(7,3)
\put(2,0.5){$r\als = r\al\mu$}
\put(8,2.7){$\Ax$}
\path(10,3)(17,3)
\put(12,0.5){$s = r\al\mu$}
\path(2,3)(2,5)
\put(1.6,2.5){$r\al\mu$}
\put(1,1.8){$=r\als$}
\path(6,3)(6,5)
\put(6,2.5){$s\mu$}
\put(6,5.2){$\Ax\als\mu$}
\path(16,3)(16,5)
\put(16,2.3){$s\mu$}
\put(16,5.2){$\Ax\als\mu$}
\path(14,4)(14,5)
\put(14,2.5){$s$}
\put(13.5,2){$=r\al\mu$}
\Thicklines
\path(3,9)(4,9)
\path(2.5,10)(3,9)
\put(2.5,8.3){$r\als$}
\path(4,9)(5,11)
\put(4,8.3){$s$}
\put(4,11.2){$\Ax\alc$}
\path(11,11)(12,9)
\put(11.2,8.3){$r\als$}
\put(10,11.2){$\Ax\alc$}
\path(12,9)(14,9)
\path(14,9)(15,11)
\put(14,8.3){$s$}
\path(2,3)(4,3)
\path(2,3)(2,4)
\put(2.2,4){$x$}
\path(2,4)(1,5)
\put(0,5.2){$\Ax\alc$}
\path(4,3)(5,5)
\put(4,2.3){$s$}
\path(11,5)(12,3)
\put(11,5.2){$\Ax\alc$}
\put(11.6,2.3){$r\als$}
\path(12,3)(14,3)
\path(14,3)(14,4)
\put(14.2,4){$y$}
\path(14,4)(13,5)
\end{picture}
\caption{$\Ax\alc \cap \Ax\ma$}
{\label{I-fig4}}
\end{figure}
Since $r$ and $s$ are not necessarily finite, we will show a
contradiction even when one of $r$ and $s$ is not finite. Recall
that at least one of $r$ and $s$ is finite.

Case (1), $[r\als, s]=[r\als\ka\alc, s\ka\alc]$.
\begin{itemize}
\item $r\als =r\als\ka\alc \Leftrightarrow r\ial= r\als\ka$. But
$r\ial \prec r \prec r\als\ka$. Contradiction.
\item $s=s\ka\alc$. Then $s \prec s\ka \Leftrightarrow s\alc
\prec_\alc s\ka\alc=s$. But $s \preceq_\al s\al \Leftrightarrow
s\ial \preceq s \Leftrightarrow s\als \preceq_\alc s\alc $. Since
$s \prec_\alc  s\als \prec_\alc s\alc$, we get a contradiction.
\end{itemize}

Case (2), $[r\al\mu,s]=[r\als\ka\alc, s\ka\alc]$.
\begin{itemize}
\item $r\al\mu=r\als\ka\alc$. Contradiction as shown in case(2) of (I).
\item $s=s\ka\alc$. Contradiction as in Case (1).
\end{itemize}

Case (3), $[x,s] = [r\als\ka\alc, s\ka\alc]$ and  $r\als=r\al\mu$.
\begin{itemize}
\item $x = r\als\ka\alc \preceq_\alc r\als \Rightarrow r\als\ka
\preceq r\ial$. But $r\ial \prec r\als \prec r\als\ka$.
Contradiction.
\item $s=s\ka\alc$. Contradiction as in Case (1).
\end{itemize}

Case (4), $[s,y] =[r\als\ka\alc, s\ka\alc]$ and $s=r\al\mu$. Then
both $r$ and $s$ are finite. We have $r\al\mu=s=r\als\ka\alc$.
Contradiction as in Case (2) of (I).

For (II), note that $\Ax \cap \Ax\alc =
\emptyset$. And every arguments in this case reduce to the case
(I) with $\Ax \cap \Ax\alc = \emptyset$ and $s \prec s\mu
=r\al\mu$ (Fig \ref{fig:II}).

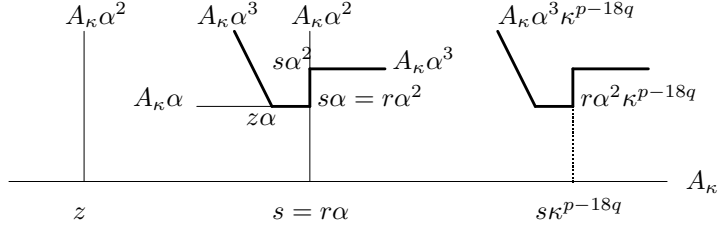
\begin{figure}[htb]
\setlength{\unitlength}{0.5cm}
\begin{picture}(20,6)
\path(0,1)(17.5,1)
\put(18,0.8){$\Ax$}
\path(2,1)(2,5)
\put(1.5,5.2){$\Ax\als$}
\put(1.7,0){$z$}
\path(8,1)(8,5)
\put(7.5,5.2){$\Ax\als$}
\dottedline{.1}(15,1)(15,3)

\path(8,3)(5,3)
\put(3.3,3){$\Ax\al$}
%\put(1.7,0.5){$r\al$}
\put(7,0){$s=r\al$}
\put(8.2,3){$s\al=r\als$}

\put(7,4){$s\als$}
\put(15.2,3){$r\als\ka$}

\put(14,0){$s\ka$}
\put(10.2,4){$\Ax\alc$}

\Thicklines
%\path(8,2)(10,2)
\path(8,4)(10,4)
\path(8,3)(8,4)
\path(8,3)(7,3)
\put(6.3,2.5){$z\al$}
\path(7,3)(6,5)

\put(5,5.2){$\Ax\alc$}

\path(15,4)(17,4)
\path(15,3)(15,4)
\path(15,3)(14,3)
\path(14,3)(13,5)
\put(13,5.2){$\Ax\alc\ka$}

\end{picture}
\caption{Case (II) $s=r\al$}{\label{fig:II}}
\end{figure}
For (III), note that $r\al \preceq_\als s\al \prec_\als r\als
\preceq_\als s\als$. In particular, $\Ax \cap \Ax\alc= \emptyset$
(Fig \ref{fig-II}).

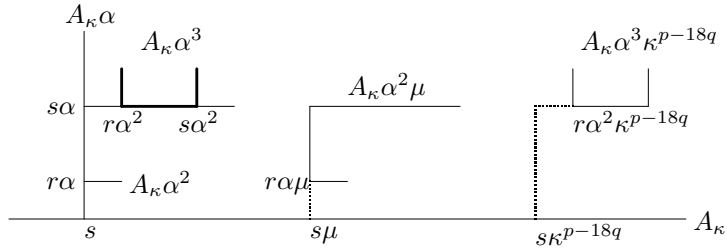
\begin{figure}
\setlength{\unitlength}{0.5cm}
\begin{picture}(20,7)

\path(0,1)(18,1)
\put(18.2,0.7){$\Ax$}
\path(2,1)(2,6)
\put(1.5,6.2){$\Ax\al$}
\put(2,0.5){$s$}
\path(2,2)(3,2)
\put(1,1.8){$r\al$}
\put(3.2,1.7){$\Ax\als$}
\path(2,4)(6,4)
\put(1,3.8){$s\al$}

\path(8,2)(8,4)
\put(6.8,1.8){$r\al\mu$}
\put(8,0.5){$s\mu$}
\dottedline{0.1}(8,1)(8,2)
\path(8,2)(9,2)
\path(8,4)(12,4)
\put(9,4.3){$\Ax\als\mu$}

\path(15,4)(17,4)
\put(15,3.3){$r\als\ka$}
\put(14,0.3){$s\ka$}
\dottedline{0.1}(14,1)(14,4)
\dottedline{0.1}(14,4)(15,4)
\path(15,4)(15,5)
\put(15.2,5.5){$\Ax\alc\ka$}
\path(17,4)(17,5)

\Thicklines
\path(3,4)(5,4)
\put(3.5,5.5){$\Ax\alc$}
\put(2.5,3.3){$r\als$}
\put(4.5,3.3){$s\als$}
\path(3,4)(3,5)
\path(5,4)(5,5)
\end{picture}
\caption{Case (III)}
{\label{fig-II}}
\end{figure}
The bridge from $\Ax\alc$ to $\Ax\ak$ is $[s\ka, r\als\ka]\alc$.
On the other hand, the bridge from $\Ax\alc$ to $\Ax\als\mu = \Ax\ma$ is
$[r\als, r\al\mu]$. So we have
\begin{itemize}
\item $s\ka\alc=r\als$ and
\item $r\al\mu=r\als\ka\alc.$
\end{itemize}
Hence $r$ and $s$ are finite and we get a contradiction as shown
in \ref{eqn:contradiction}.

Now we may assume the map $\al : (\Ax, \prec) \rightarrow (\Ax\al,
\prec_\al)$ is order reversing. Equivalently, $s\al \prec_\al
r\al$.

If $[r,s] \cap [s\al,r\al]\ne \emptyset$, intermediated value
theorem is applied to show $Nonsep(\al) \cap \Ax \ne \emptyset$.
Hence we assume $[r,s] \cap [s\al, r\al] =\emptyset$. By symmetry,
we can assume
$$ r \prec_\al s \prec_\al s\al \prec_\al r\al \text{ (and hence $s$ is finite)}.$$
\begin{figure}[hbt]
\setlength{\unitlength}{0.5cm}
\begin{picture}(20,7)
\path(0,1)(19,1)
\put(19.2,1){$\Ax$}
\path(5,1)(7,6)
\put(6.5,6.2){$\Ax\al$}
\put(4.8,0.5){$s$}
\path(5.4,2)(1,2)
\put(0,2){$\Ax\als$}
\put(5.6,2){$s\al$}
\path(6.6,5)(4,6)
\put(3,6){$\Ax\als$}
\put(6.8,5){$r\al$}

\path(12.4,2)(13.6,5)
\path(13.6,5)(11,6)
\put(10,6){$\Ax\als\mu$}
\path(12.4,2)(8,2)
\put(12.6,2){$s\al\mu$}
%\path(6,1)(10.6,1)
\put(11.8,0.5){$s\mu$}
\dottedline{0.1}(12,1)(12.4,2)
\dottedline{0.1}(17,2)(18.4,2)
\dottedline{0.1}(18.4,2)(18,1)
\path(15,2)(17,2)
\put(15.2,4){$\Ax\alc\ka$}
\put(17.2,2.1){$s\als\ka$}
\path(15,2)(15,4)
\path(17,2)(17,4)
\put(17,0.3){$s\ka$}

\Thicklines
\path(2,2)(4,2)
\put(1.5,1.3){$r\als$}
\put(3.5,1.3){$s\als$}
\path(2,2)(2,4)
\path(4,4)(4,2)
\put(2.5,4){$\Ax\alc$}
\end{picture}
\caption{order reversing}{\label{order reversing}}
\end{figure}
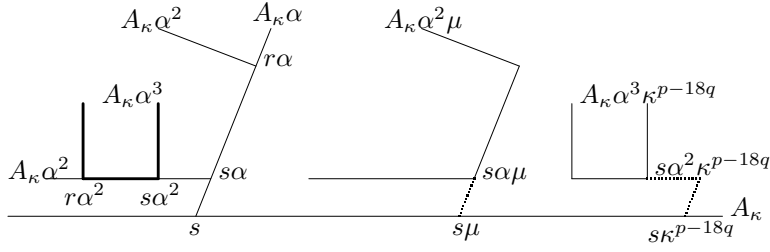

The bridge from $\Ax\alc$ to $\Ax\ak$ is $[s\ka\alc,
s\als\ka\alc]$(Fig \ref{order reversing}).
The bridge from $\Ax\alc$ to $\Ax\als\mu$ is
$[s\als, s\al\mu]$. So we have
\begin{align*}
&s\als\ka\alc = s\al\mu\\
&\Leftrightarrow s\be\als = s\ibes\\
&\Leftrightarrow s\be\als\bes=s\\
&\Leftrightarrow s\ial\mu\al=s \quad\text{by (\ref{eqn:rel})}\\
&\Leftrightarrow s\ial\mu = s\ial.
\end{align*}
 But $s\ial \prec s\ial\mu$. Contradiction.
\end{pf}

The line of reasoning used in this section shows that one actually
has
\begin{lem}{\label{lem:lem7.4}}
Suppose $Y$ is a $\kappa$-invariant embedded copy of $\Real$ in
$\mathcal{T}$ on which $\kappa$ acts freely. If
\begin{itemize}
\item $\emptyset \ne Y \cap Y\al \subset [r,s]$ for some $r,s \in
Y$, or
\item $Y \cap Y\al = \emptyset$ and the bridge from $Y$ to $Y\al$
has the form $[[r,s]]$ for some $r\sim r' \in Y$, $s \sim s' \in
Y\al$,
\end{itemize}
then the action $\G$ on $\mathcal{T}$ has a global fixed point.
\end{lem}
%-----------------------------------------------------------------

\section{$ Nonsep(\kappa) \ne \emptyset$}{\label{sec:nonempty}}
We will complete the proof of Main Theorem by showing below that
there is no nontrivial action on $\mathcal{T}$ when
$Nonsep(\kappa)\ne \emptyset$

\begin{lem}{\label{lem:lem8.1}}
There is no $x \in \mathcal{T}$ which is nonseparated by $\kappa$
and $\al$.
\end{lem}

\begin{lem}{\label{lem:ka}}
%Suppose that $Nonsep(\al)=\emptyset$.
If $Fix(\ka) \cap Nonsep(\kappa) \cap C_\al \ne \emptyset$, then
the action $\G$ on $\mathcal{T}$ is trivial.
\end{lem}

\begin{pf}
Let $x \in Fix(\ka) \cap Nonsep(\kappa) \cap C_\al$. By Lemma
\ref{lem:4.6}, $x$ and $x\al$ are comparable. If $x\al=x$,
$x=x\ka=x\ials\be\ial=x\be\ial \Leftrightarrow x=x\be$. Because
$\G =<\al, \be>$, instead we can assume that $x < x\al$. Then
\begin{align*}
\begin{split}
&x = x \ka = x \ials \be \ial < x\be \ial\\
&\Rightarrow x\al < x\be\\
&\Rightarrow x < x\be\\
&\Rightarrow x < x\bes = x \mu^{-1}\ial\\
&\Rightarrow x\al < x\mu^{-1}\\
&\Rightarrow x < x\mu^{-1}.
\end{split}
\end{align*}

Since $x \sim x \mu$, we get a contradiction.
\end{pf}
\begin{lem}
If $Nonsep(\kappa)\cap C_\al \ne\emptyset$, then the action $\G$
on $\mathcal{T}$ is trivial.
\end{lem}

\begin{pf}
Let $x \in Nonsep(\kappa)\cap C_\al$. By Lemma \ref{poset}, we may
assume that $x \sim x\kappa$ but $x \ne x\kappa$. Set
$\mathcal{T}_0$ = $\mathcal{T}_{\{x,x\kappa\}}$. We assume that $x
< x\al$.

If $x\ial \in \mathcal{T}_0$ or $x\al \in \mathcal{T}_0$, then the
ideal point $\hat{x} \in \widehat{\mathcal{T}}$ is fixed by
$\kappa$ and related to $\hat{x}\al$, and Lemma \ref{poset}
applies.

So we may assume that $x\al, x\ial \not\in \mathcal{T}_0$. Since
$x < x\al$, either  $\mathcal{T}_0 \subset x^+$ and $x\al \in y^-$
for some $y \sim x, y \ne x$ or $\mathcal{T}_0 \subset x^-$ and
$x\al \in y^+$ for some $y \sim x, y \ne x$. We may assume the
first possibility holds. Note that $\{x,y\} \subset [[x\ial,
x\al]]$ and so $d(x,x\al)=2n > 0$.  In particular, we have
$Nonsep(\al) =\emptyset$ and $C_\al = \Aal$ by Lemma
\ref{lem:cor4.11} (Fig \ref{fig:C}). We can also assume that $x
\ne x\ka$ and $y \ne y\ka$ by Lemma \ref{lem:ka}.

\begin{figure}[hbt]
\setlength{\unitlength}{0.5cm}
\begin{picture}(15,2)
\put(2,1){\circle*{0.2}}
\put(0.7,0){$x\ial$}
\put(3,1){\circle*{0.2}}
\put(3,0){$y\ial$}
\put(7,1){\circle*{0.2}}
\put(6.5,0){$x$}
\put(8,1){\circle*{0.2}}
\put(8,0){$y$}
\put(12,1){\circle*{0.2}}
\put(11.5,0){$x\al$}
\put(13,1){\circle*{0.2}}
\put(13,0){$y\al$}

\put(0,1){\vector(1,0){1.9}}
\put(5,1){\vector(-1,0){1.9}}
\put(5,1){\vector(1,0){1.9}}
\put(10,1){\vector(-1,0){1.9}}
\put(10,1){\vector(1,0){1.9}}
\put(15,1){\vector(-1,0){1.9}}
\end{picture}
\caption{$C_\al$}{\label{fig:C}}
\end{figure}
Then we have 3 cases.
\begin{enumerate}
\item $y\mu=y$ and $x \ne y\ka$.

$d(y, y\ma) = d(y,y\als)=4n$. But $d(y,y\ak) = d(y,y\ka\alc) +
d(y\ka\alc,y\ak) = 6n+6n = 12n$ (Fig \ref{fig:case1-1}). Since
$n>0$, this is impossible.

\begin{figure}[hbt]
\setlength{\unitlength}{0.5cm}
\begin{picture}(20,7)
\put(2,1){\circle*{0.2}}
\put(2,0){$x$}
\put(3,1){\circle*{0.2}}
\put(3,0){$y$}
\put(9,1){\circle*{0.2}}
\put(8,0){$x\alc$}
\put(10,1){\circle*{0.2}}
\put(10,0){$y\alc$}
\put(2.5,2){\circle*{0.2}}
\put(3,2){$y\ka$}
\put(9.5,2){\circle*{0.2}}
\put(10,2){$y\ka\alc$}
\put(2.5,4){\circle*{0.2}}
\put(3,3.6){$x\alc\ka$}
\put(2.5,5){\circle*{0.2}}
\put(3,5){$y\alc\ka$}
\put(9.5,4){\circle*{0.2}}
\put(10,3.6){$x\ak$}
\put(9.5,5){\circle*{0.2}}
\put(10,5){$y\ak$}

\put(0,1){\vector(1,0){1.9}}
\put(6,1){\vector(-1,0){2.9}}
\put(6,1){\vector(1,0){2.9}}
\put(15,1){\vector(-1,0){4.9}}
\put(2.5,3){\vector(0,-1){0.9}}
\put(2.5,3){\vector(0,1){0.9}}
\put(9.5,3){\vector(0,-1){0.9}}
\put(9.5,3){\vector(0,1){0.9}}
\put(2.5,7){\vector(0,-1){1.9}}
\put(9.5,7){\vector(0,-1){1.9}}
\end{picture}
\caption{$x \ne y\ka$}{\label{fig:case1-1}}
\end{figure}

\item $y\mu=y$ and $x = y\ka$.
    \begin{enumerate}
    \item $y\be \in \Aal \Leftrightarrow y\als\ka\in\Aal$.

    Then $d(x, y\als\ka)=d(y,y\als) = 4n$. Since $y\als\in y^-$, $y\als\ka \in (y\ka)^-=
    x^-$. Moreover $y\als\ka \in \Aal$. Therefore,
    \begin{align*}
    \begin{split}
    &y\als\ka = x\ials\\
    &\Leftrightarrow y\be = x\ial\\
    &\Rightarrow y\be < x\\
    &\Leftrightarrow y\ial = y\mu^{-1}\ial =y\bes < x\be.
    \end{split}
    \end{align*}
    Since $x\be \sim y\be=x\ial \sim y\ial$, we get a contradiction.

    \item $y\be \not\in\Aal$.

    Since $d(x\al, y\be)= d(x\al, y\als\ka\al)=d(x,y\als\ka)=d(y, y\als)
    = 4n = d(x\al,x\ial)$, $x\ial \not\in [[x\al, y\be]]$.

    If $x\in[[x\al, y\be]]$, then $y\be <x$ (Fig \ref{fig:1-2}).
    Therefore $x\be>y\ial=y\imu\ial=y\bes $. Contradiction.

    If $x\not\in[[x\al,x\be]]$,
    $y>y\be \Rightarrow y > y\bes=y\ial$(Fig \ref{fig:1-2}).
    Contradiction.

    \begin{figure}[htb]
\setlength{\unitlength}{0.5cm}
\begin{picture}(20,4)
\put(1,3){\circle*{0.2}}
\put(0,3.5){$x\ial$}
\put(2,3){\circle*{0.2}}
\put(2,3.5){$y\ial$}
\put(4,3){\circle*{0.2}}
\put(3.8,3.5){$x$}
\put(5,3){\circle*{0.2}}
\put(5,3.5){$y$}
\put(7,3){\circle*{0.2}}
\put(6.2,3.5){$x\al$}
\put(8,3){\circle*{0.2}}
\put(8,3.5){$y\al$}
\put(1,1){\circle*{0.2}}
\put(1.5,1){$x\be$}
\put(0.5,0.5){\circle*{0.2}}
\put(0.5,0){$y\be$}

\put(0,3){\vector(1,0){0.9}}
\put(3,3){\vector(-1,0){0.9}}
\put(3,3){\vector(1,0){0.9}}
\put(6,3){\vector(-1,0){0.9}}
\put(6,3){\vector(1,0){0.9}}
\put(9,3){\vector(-1,0){0.9}}
\put(3,3){\vector(-1,-1){2}}
\put(0,0){\vector(1,1){0.5}}

\path(9.5,4)(9.5,0)
\put(11,3){\circle*{0.2}}
\put(10,3.5){$x\ial$}
\put(12,3){\circle*{0.2}}
\put(12,3.5){$y\ial$}
\put(14,3){\circle*{0.2}}
\put(13.8,3.5){$x$}
\put(15,3){\circle*{0.2}}
\put(15,3.5){$y$}
\put(17,3){\circle*{0.2}}
\put(16.2,3.5){$x\al$}
\put(18,3){\circle*{0.2}}
\put(18,3.5){$y\al$}
\put(14,1){\circle*{0.2}}
\put(14.5,1){$x\be$}
\put(13.5,0.5){\circle*{0.2}}
\put(13.5,0){$y\be$}

\put(10,3){\vector(1,0){0.9}}
\put(13,3){\vector(-1,0){0.9}}
\put(13,3){\vector(1,0){0.9}}
\put(16,3){\vector(-1,0){0.9}}
\put(16,3){\vector(1,0){0.9}}
\put(19,3){\vector(-1,0){0.9}}
\put(16,3){\vector(-1,-1){2}}
\put(13,0){\vector(1,1){0.5}}
\end{picture}
\caption{$x\be$ and $\Aal$}{\label{fig:1-2}}
\end{figure}
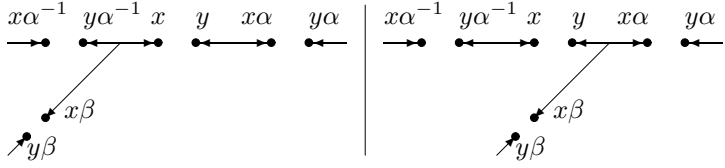

    \end{enumerate}
\item $y \ne y\mu$.

Since $x\mu\als \sim y\als \in y^-$, $x\mu\als \in y^-$. Thus
$x\ma \in (y\mu)^-\subset y^+$. But $x\ak \in y^-\cup\{y\}$(see
Fig \ref{fig:case1-1}).

\end{enumerate}
\end{pf}

\begin{lem}{\label{lem:lem8.5}}
If $Nonsep(\kappa) \ne \emptyset$ and $Nonsep(\al) \cap C_\kappa
\ne \emptyset$, then the action $\G$ on $\mathcal{T}$ is trivial.
\end{lem}
\begin{pf}
Let $x \in Nonsep(\al)\cap C_\kappa$. By Lemma \ref{lem:lem4.9},
either $x\in Fix(\al)$ or $x$ lies on some local axis $\Ax^i
\approx \Real$ (in $\mathcal{T}$) for $\kappa$. By Lemma
\ref{lem:lem8.1}, we may assume that $x$ lies on some local axis
$\Ax^i$. Then either $x\in Fix(\al)$ or the ideal point $\hat{x}
\in \widehat{\mathcal{T}}$ is fixed by $\al$ and related to
$\hat{x}\kappa$. In either case, Lemma \ref{poset} applies.
\end{pf}

\begin{lem}{\label{lem:lem8.6}}
If $G(p,q)$ acts nontrivially on $\mathcal{T}$, then:
\begin{itemize}
\item $C_\kappa \cup Nonsep(\kappa) \subset X_{j_0}$ for some $j_0
\in \mathcal{J}$, and
\item $C_\al \cup Nonsep(\al) \subset T_{i_0}$ for some $i_0
\in \mathcal{I}$.
\end{itemize}
\end{lem}

\begin{pf}
By Lemma \ref{lem:lem8.5}, $(C_\kappa \cup Nonsep(\kappa)) \cap
Nonsep(\al) = \emptyset$. By Lemma \ref{lem:cor4.12} therefore,
$C_\kappa \cup Nonsep(\kappa) \subset X_{j_0}$ for some $j_0 \in
\mathcal{J}$. A symmetric argument proves the second statement.
\end{pf}

\begin{prop}
Suppose $Nonsep(\kappa) \ne \emptyset$. Then the action is
trivial.
\end{prop}
\begin{pf}
Let $i_0, j_0$ be as guaranteed in Lemma \ref{lem:lem8.6}. Suppose
first that $T_{i_0} \kappa =T_{i_0}$. As remarked above,
$\Ax^{i_0} \approx \Real$. By Lemma \ref{lem:lem8.6},
$Nonsep(\kappa) \cup \Ax^{i_0} \subset X_{j_0}$.

Consider first the possibility that $X_{j_0}\al = X_{j_0}$, and
hence $\Aal^{j_0} \subset T_{i_0}$. In fact, $T_{i_0} \cap
X_{j_0}$ is a subtree of $\mathcal{T}$ containing both $\Ax^{i_0}$
and $\Aal^{j_0}$. Therefore, if $\Ax^{i_0} \cap
\Aal^{j_0}=\emptyset$, the bridge from $\Ax^{i_0}$ to $\Aal^{j_0}$
lies in $T_{i_0} \cap X_{j_0}$. If either of the two potential
endpoints of $\Ax^{i_0}$(respectively, $\Aal^{j_0}$) exist in
$\mathcal{T}$, they are in $Nonsep(\kappa)$ (respectively,
$Nonsep(\al)$) and hence are not elements of $\mathcal{T}_{i_0}$
(respectively, $X_{j_0}$), and therefore cannot be on the bridge.
Hence this bridge has the form $[[u,v]]$ or $[[u,v))$, where $u$
and $v$ are not separated from points in $\Ax^{i_0}$ and
$\Aal^{j_0}$ respectively. Computing $\Ax^{i_0}$ in this case, we
see that $\Ax^{i_0} \cap \Ax^{i_0}\al=\emptyset$, with the bridge
from $\Ax^{i_0}$ to $\Ax^{i_0}\al$ of the form $[[u,w]]$ for some
$w \sim w'\in\Ax^{i_0}\al$. So Lemma \ref{lem:lem7.4} reveals that
the action of $G(p,q)$ on $\mathcal{T}$ is necessarily trivial. On
the other hand, if $\Ax^{i_0} \cap \Aal^{j_0} \ne \emptyset$ then
Lemma \ref{lem:lem8.1} guarantees that $ \Ax^{i_0} \cap
\Aal^{j_0}\subset [u,v]$ for some $u,v \in \Ax^{i_0}$. Computing
$\Ax^{i_0}\al$ in this case, we see that one of the two conditions
of Lemma \ref{lem:lem7.4} is satisfied, and so once again, the
action of $G(p,q)$ on $\mathcal{T}$ must be trivial.

Next consider the possibility that $X_{j_0}\al =X_{j_1}\ne
X_{j_0}$. Let $y$ and $y\al$ denote the roots of $X_{j_0}$ and
$X_{j_0}\al$, respectively. Let $[[y,r]]$ denote the bridge from
$y$ to $\Ax^{i_0}$ in $\mathcal{T}$. By Lemma \ref{lem:lem8.1}, we
may assume that $r \sim r'$ for some $r'\in\Ax^{i_0}$. So
$\Ax^{i_0} \cap \Ax^{i_0}\al=\emptyset$ with bridge $[[r, r\al]]$.
Again, by Lemma \ref{lem:lem7.4}, the action of $G(p,q)$ on
$\mathcal{T}$ has a global fixed point.

 Suppose that $T_{i_0}\kappa = T_{i_1} \ne T_{i_0}$.
In particular, $x \ne x\kappa$. As shown in the proof of
Proposition 8.7 in \cite{RSS}, we have $x\al \in x^-$ and $x\in
(x\al)^-$ (Fig \ref{fig:x}).

\begin{figure}[htb]
\setlength{\unitlength}{0.5cm}
\begin{picture}(10,2)
\put(1,1){\circle*{0.2}}
\put(9,1){\circle*{0.2}}
\put(1,0){$x\al$}
\put(9,0){$x$}
\put(5,1){\vector(-1,0){3.9}}
\put(5,1){\vector(1,0){3.9}}
\end{picture}
\caption{$x$ and $x\al$ when $T_{i_0}\kappa =
T_{i_1}$}{\label{fig:x}}
\end{figure}
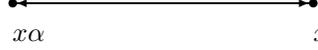
If $x=x\mu =x\ka$, then $x=x\kappa^p$. Since $(p,q)=1$, $x=x\kappa$.
Hence $x$ is not equal to at least one of $x\mu$ and $x\ka$.

Suppose that $x=x\als$. Assume first $x=x\ka$ (and hence $x\ne
x\mu$). Then $x\al\imu=x\alc\imu=x\ialc\mu\als=x\al\mu\als.$ Since
$d(x\imu, x\al\imu) =d(x,x\al) =d(x\mu,x\al\mu)=d(x\mu\als,
x\al\mu\als)$ and $x\imu \sim x =x\als \sim x\mu\als$,
\begin{align*}
\begin{split}
&x\imu=x\mu\als \mbox{ (see Fig \ref{fig:als=ka})}\\
&\Leftrightarrow x\als=x = x\mu\als\mu = x\al\be\als\\
&\Leftrightarrow x = x\al\be\\
&\Leftrightarrow x\be = x\als\ka\al = x \al = x\ibe \\
&\Rightarrow x\imu\ial = x\bes =x\\
&\Leftrightarrow x\imu = x\al.
\end{split}
\end{align*}
Since $x\imu \sim x$, we get a contradiction.
\begin{figure}
\setlength{\unitlength}{0.5cm}
\begin{picture}(15,2)
\put(2,1){\circle*{0.2}}
\put(2,0){$x\al$}
\put(4.5,1){\vector(-1,0){2.4}}
\put(4.5,1){\vector(1,0){2.4}}
\put(7,1){\circle*{0.2}}
\put(5,0){$x\als=x$}
\put(8,1){\circle*{0.2}}
\put(8,0){$x\imu$}
\put(10.5,1){\vector(-1,0){2.4}}
\put(10.5,1){\vector(1,0){2.4}}
\put(13,1){\circle*{0.2}}
\put(12,0){$x\al\imu=x\al\mu\als$}

\end{picture}
\caption{$x=x\als=x\ka$}{\label{fig:als=ka}}
\end{figure}
Now we can assume

Note that
$x\ak\in x^-$ (Fig \ref{fig:case1}). But $x\ma \sim x\als\mu=x\mu\sim x$.
So $x\ma \in x^+ \cup \{x\}$. Contradiction.

Now we may assume $x\ne x\als$. Suppose that $x=x\alc$. Then
$$x\ak=x\ka\alc\sim x\alc=x.$$
But $x\ma \not\sim x$ as shown in Fig \ref{fig:case1}(when $x\ne
x\mu$) and Fig \ref{fig:case3}(when $x=x\mu$). Hence we also
assume that $x\ne x\alc$.

Recall that $x$ is not equal to at
least one of $x\mu$ and $x\ka$.

We have  3 cases. Note that
$x\als, x\alc \in x^-$ and $x \in (x\als)^- \cap (x\alc)^-$.
\begin{enumerate}
\item $x \ne x\ka$ and $x \ne x\mu$. \\
%$x\al \in x^-\Rightarrow x\als, x\alc \in x^-$.
 As shown in Fig \ref{fig:case1},  $x\ma \in x^+$ but
$x\ak \in x^-$.
\item $x \ne x\mu$ and $x = x\ka$\\
It follows that $x\alc\mu^{-1} = x\al^{-3}\mu\als$. Then
$x\alc\mu^{-1} \in x^+$, but $x\al^{-3}\mu\als \in x^-$ (see Fig
\ref{fig:case2}).
\item $x \ne x\ka$ and $x=x\mu$.\\
As shown in the figure of case (1), $x\ak \in (x\alc)^+$. But
$x\ma=x\als\mu \in (x\alc)^-$ (Fig \ref{fig:case3}).
%\item $x=x\mu$ and $x=x\ka$.\\
%Then $x=x\kappa^p$. Since $(p,q)=1$, $x=x\kappa$. But we are in
%case $T_{i_0}\kappa = T_{i_1}\ne T_{i_0}$.
\end{enumerate}

\begin{figure}
\setlength{\unitlength}{0.5cm}
\begin{picture}(20,6)
\put(2,1){\circle*{0.2}}
\put(0,0){$x\ak$}
\put(6,1){\circle*{0.2}}
\put(3.6,1.2){$x\ka\alc$}
\put(4,1){\vector(-1,0){1.9}}
\put(4,1){\vector(1,0){1.9}}
\put(7,1){\circle*{0.2}}
\put(7,0){$x\alc$}
\put(11,1){\circle*{0.2}}
\put(10.8,0){$x$}
\put(9,1){\vector(-1,0){1.9}}
\put(9,1){\vector(1,0){1.9}}
\put(12,1){\circle*{0.2}}
\put(12,0){$x\ka$}
\put(16,1){\circle*{0.2}}
\put(15,1.2){$x\alc\ka$}
\put(14,1){\vector(-1,0){1.9}}
\put(14,1){\vector(1,0){1.9}}

\put(2,5){\circle*{0.2}}
\put(1,4){$x\mu\als$}
\put(1,5){\vector(1,0){0.9}}
\put(3,5){\circle*{0.2}}
\put(3,4){$x\als$}
\put(7,5){\circle*{0.2}}
\put(6.5,4){$x$}
\put(5,5){\vector(-1,0){1.9}}
\put(5,5){\vector(1,0){1.9}}
\put(8,5){\circle*{0.2}}
\put(8,4){$x\mu$}
\put(12,5){\circle*{0.2}}
\put(10.5,4){$x\als\mu$}
\put(10,5){\vector(-1,0){1.9}}
\put(10,5){\vector(1,0){1.9}}
\put(13,5){\circle*{0.2}}
\put(13,4){$x\ma$}
\put(15,5){\vector(-1,0){1.9}}

\end{picture}
\caption{Case (1) $x\ne x\ka, x\ne x\mu$}{\label{fig:case1}}
\end{figure}
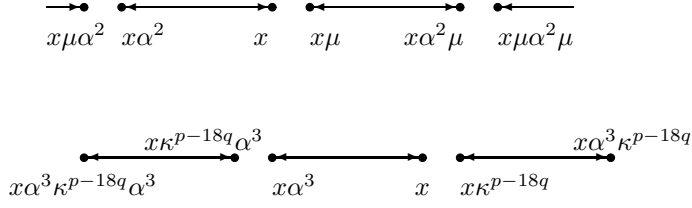

\begin{figure}
\setlength{\unitlength}{0.5cm}
\begin{picture}(20,8)
\put(1,3){\circle*{0.2}}
\put(0,2){$x\al^{-3}\mu\als$}
\put(2,3){\vector(-1,0){0.9}}
\dashline[+30]{0.2}(2,3)(7,3)
\put(7,3){\circle*{0.2}}
\put(6,2){$x\mu\als$}
\put(4,3){\vector(1,0){2.9}}
\put(8,3){\circle*{0.2}}
\put(8,2){$x\als$}
\put(13,3){\circle*{0.2}}
\put(12.5,2){$x$}
\put(10.5,3){\vector(-1,0){2.4}}
\put(10.5,3){\vector(1,0){2.4}}
\put(9,0.5){\circle*{0.2}}
\put(9.5,0){$x\al^{-3}$}
\put(10,1.5){\vector(-1,-1){0.9}}
\dashline[+30]{0.2}(10,1.5)(11,2.5)
\put(14,3){\circle*{0.2}}
\put(14,2){$x\mu$}
\put(19,3){\circle*{0.2}}
\put(18,2){$x\al^{-3}\mu$}
\put(16.5,3){\vector(-1,0){2.4}}
\put(16.5,3){\vector(1,0){2.4}}

\put(3,7){\circle*{0.2}}
\put(2,6){$x\alc$}
\put(9,7){\circle*{0.2}}
\put(8.5,6){$x$}
\put(10,7){\circle*{0.2}}
\put(10,6){$x\mu^{-1}$}
\put(16,7){\circle*{0.2}}
\put(15,6){$x\alc\mu^{-1}$}
\put(6,7){\vector(-1,0){2.9}}
\put(6,7){\vector(1,0){2.9}}
\put(13,7){\vector(-1,0){2.9}}
\put(13,7){\vector(1,0){2.9}}
\end{picture}
\caption{Case (2) $x\ne x\mu, x=x\ka$ }{\label{fig:case2}}
\end{figure}
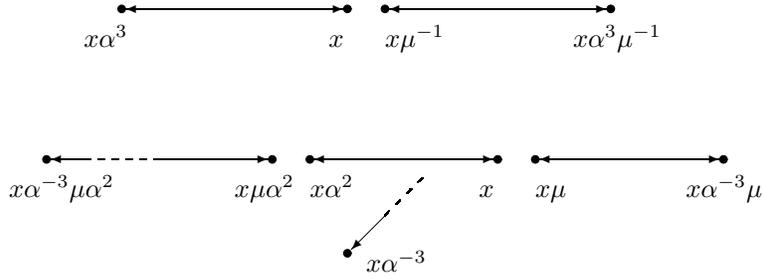

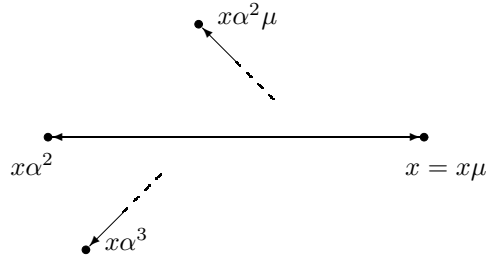
\begin{figure}[htb]
\setlength{\unitlength}{0.5cm}
\begin{picture}(15,8)
\put(2,4){\circle*{0.2}}
\put(1,3){$x\als$}
\put(12,4){\circle*{0.2}}
\put(11.5,3){$x=x\mu$}
\put(7,4){\vector(-1,0){4.9}}
\put(7,4){\vector(1,0){4.9}}

\put(3,1){\circle*{0.2}}
\put(3.5,1){$x\alc$}
\put(4,2){\vector(-1,-1){0.9}}
\dashline[30]{0.2}(4,2)(5,3)
\put(6,7){\circle*{0.2}}
\put(6.5,7){$x\als\mu$}
\put(7,6){\vector(-1,1){0.9}}
\dashline[30]{0.2}(7,6)(8,5)
\end{picture}
\caption{Case (3) $x\ne x\ka, x=x\mu$}{\label{fig:case3}}
\end{figure}
\end{pf}
% ----------------------------------------------------------------

\end{document}